\documentclass[preprint,review,12pt]{elsarticle}

\usepackage{amssymb}
\usepackage{amsmath}
\usepackage{amsthm}
\usepackage{nomencl}
\usepackage{bbm}
\usepackage{algorithm}
\usepackage{algpseudocode}
\usepackage{booktabs}
\usepackage{tabularx}
\usepackage{multirow}
\usepackage{hyperref}
\usepackage{xurl}

\newtheorem{remark}{Remark}
\usepackage{xcolor}
\usepackage{longtable}

\makenomenclature

\journal{Acta Astronautica}

\begin{document}

\begin{frontmatter}



\title{Spare Strategy for Large-Scale Satellite Constellations \\
Under Dual Resupply Channels Using Markov Chain}


\author[gti]{Seungyeop Han\footnote{Ph.D. Candidate, Daniel Guggenheim School of Aerospace Engineering}}
\author[mec]{Shoji Yoshikawa\footnote{Chief Expert, Advanced Technology R\&D Center, Mitsubishi Electric Corporation}}
\author[mec]{Takumi Noro\footnote{Researcher, Advanced Technology R\&D Center, Mitsubishi Electric Corporation}}
\author[mec]{Takumi Suda\footnote{Head Researcher, Advanced Technology R\&D Center, Mitsubishi Electric Corporation}}
\author[gti]{Koki Ho\footnote{Dutton-Ducoffe Professor, Associate Professor, Daniel Guggenheim School of Aerospace Engineering, AIAA Senior Member, kokiho@gatech.edu (Corresponding Author)}}

\affiliation[gti]{
    organization={Georgia Institute of Technology},
    addressline={Daniel Guggenheim School of Aerospace Engineering},
    city={Atlanta},
    postcode={30332},
    state={Georgia},
    country={USA}
}

\affiliation[mec]{
    organization={Mitsubishi Electric Corporation},
    addressline={Advanced Technology R\&D Center},
    city={Amagasaki},
    postcode={661-8861},
    country={Japan}
}

\begin{abstract}
\textcolor{black}{This paper presents a Markov-chain-based method for the early-phase analysis and design of hybrid spare-management architectures for large-scale satellite constellations.} The hybrid strategy combines two channels: an indirect path that stages spares in parking orbits via heavy launch for later transfer to constellation planes, and a direct path that delivers spares to in-plane orbits using small launch vehicles. \textcolor{black}{To assess the long-run viability of such concepts of operations, satellite failure and replenishment processes are modeled as a Markov chain:} the indirect channel follows a periodic-review reorder-point/order-quantity policy, while the direct channel uses a standard reorder-point/order-quantity policy. These coupled chains yield a periodic steady state over the right ascension of the ascending node cycle via fixed-point iteration, and the stationary distributions provide rigorous cost and resilience metrics. By directly modeling the stochastic, multi-echelon dynamics governed by orbital mechanics, our framework avoids the aggregation assumptions of prior works and remains valid across a wider operating domain. We also introduce an approximate analysis that preserves delay statistics while significantly reducing model size. Building on this fast, accurate analysis, we formulate a cost minimization problem with resilience constraints and solve it using a genetic algorithm. The framework is channel-neutral; the optimization autonomously selects the preferred path and roles. \textcolor{black}{A case study validates the analysis against Monte Carlo simulations and demonstrates the practical value of the framework in identifying the conditions under which the hybrid policy outperforms pure strategies.}
\end{abstract}




\end{frontmatter}


\nomenclature[a01]{$\tau_\text{mc}$}{Time step of the discrete-time Markov process, in days}
\nomenclature[a02]{$\tau_\text{c}$}{RAAN alignment period from the perspective of in-plane (constellation) orbits}
\nomenclature[a03]{$\tau_\text{p}$}{RAAN alignment period from the perspective of parking orbits}
\nomenclature[a04]{$\tau_{\text{lv,d}}$}{Launch order processing time of the direct channel (using a small launch vehicle), in days}
\nomenclature[a05]{$\tau_{\text{lv,i}}$}{Launch order processing time of the indirect channel (using a heavy launch vehicle), in days}
\nomenclature[a06]{$\mu_{\text{lv,d}}$}{Mean interval between launches in the direct channel, in days}
\nomenclature[a07]{$\mu_{\text{lv,i}}$}{Mean interval between launches in the indirect channel, in days}
\nomenclature[a08]{$\lambda_\text{sat}$}{Failure rate of a satellite, in failures per unit time}

\nomenclature[a09]{$q_{\text{c,d}}$}{Replenishment quantity for in-plane spares via the direct channel}
\nomenclature[a10]{$q_{\text{c,i}}$}{Replenishment quantity for in-plane spares via the indirect channel}
\nomenclature[a11]{$q_{\text{p}}$}{Replenishment quantity for parking-orbit stock, in batches}
\nomenclature[a12]{$r_{\text{c,d}}$}{Reorder point for in-plane spares in the direct channel}
\nomenclature[a13]{$r_{\text{c,i}}$}{Reorder point for in-plane spares in the indirect channel}
\nomenclature[a14]{$r_{\text{p}}$}{Reorder point for parking spares, in batches}

\nomenclature[a15]{$\bar{N}_{\text{sat}}$}{Nominal number of operational satellites per in-plane orbit}
\nomenclature[a16]{$N_{\text{sat}_\text{c}}$}{Maximum in-plane state level, including operational and spare satellites}
\nomenclature[a17]{$N_{\text{sat}_\text{p}}$}{Maximum parking-orbit state level}
\nomenclature[a18]{$N_{\text{orbit}_{\text{c}}}$}{Number of in-plane orbital planes}
\nomenclature[a19]{$N_{\text{orbit}_{\text{p}}}$}{Number of parking orbital planes}

\nomenclature[a20]{$P_f$}{Failure transition matrix for in-plane states}
\nomenclature[a21]{$P_{f_{\text{p}}}$}{Demand-induced failure transition matrix for parking states}
\nomenclature[a22]{$P_{q_{\text{c,d}}}$}{Replenishment transition matrix for in-plane spares via the direct channel}
\nomenclature[a23]{$P_{q_{\text{c,i}}}$}{Replenishment transition matrix for in-plane spares via the indirect channel}
\nomenclature[a24]{$P_{q_{\text{p}}}$}{Replenishment transition matrix for parking spares}

\nomenclature[a30]{$\tilde{\pi}/\pi$}{Relative and full state probability distribution vector, $\pi = \tilde{\pi}/(\mathbbm{1}\cdot \tilde{\pi})$}
\nomenclature[a31]{$\pi^{\text{io}_{\text{c}}}/\pi^{\text{io}_{\text{p}}}$}{Expected in-plane/parking state distribution during the inter-order (IO) period}
\nomenclature[a32]{$\pi^{\text{lt}_{\text{c}}}/\pi^{\text{lt}_{\text{p}}}$}{Expected in-plane/parking state distribution during the lead-time (LT) period}
\nomenclature[a33]{$\pi^{\text{rc}_{\text{c}}}/\pi^{\text{rc}_{\text{p}}}$}{Expected in-plane/parking state distribution over the full replenishment cycle (RC)}

\printnomenclature


\section{Introduction} \label{sec1}
Compared with other distributed satellite systems, operating very large
constellations in Low Earth Orbit (LEO) is a logistics problem because of their
scale~\cite{selva2017distributed}. Large fleets face more failures, and low-cost designs may have limited redundancy~\cite{BOUWMEESTER2022108288}. To maintain service,
operators can replace failed satellites with spares~\cite{cornara1999satellite}
or attempt on-orbit servicing~\cite{luu2022orbit}. For LEO systems, replacement is often more practical because satellites are produced in large numbers, whereas making small satellites serviceable adds design and operational complexity. {The resulting
design problem is to select spare-management architecture parameters, such as
staging locations, replenishment channels, reorder thresholds, batch sizes, and
parking-orbit configuration. Such choices define the spare-management
architecture and motivate systematic early-phase analysis.}

This spare-management design problem becomes more important as LEO constellations scale. Small satellites increasingly rely on compact electronics and commercial off-the-shelf parts. At the same time, lower launch prices make frequent replenishment more feasible~\cite{potter2023mooreslaw}. This makes it 
feasible to deploy thousands of small, low-cost satellites for global services.
Examples include OneWeb, Starlink, Project Kuiper, and
Guowang~\cite{oneweb,starlink1,starlink2,kuiper2025,guowang2025}. These systems
rely on advanced constellation architecture, frequency coordination, and
inter-satellite communication. This growth reflects demand for continuous
coverage and high-capacity networks; see~\cite{Pachler2021Comp,Pachler2024Comp}
for comparative analyses of their designs.

\textcolor{black}{Constellation sustainment problems have been tackled using probabilistic models. Representative works include Ref.~\cite{1966dishon}, which used $(s,S)$ inventory modeling formulation with exponential lead times, where $s$ is the reorder point and $S$ is the order-up-to level, and Ref.~\cite{sumter2003optimal}, which used discrete-time Markov decision processes for selecting cost-effective replacement actions over finite horizons. For mega-constellations, Jakob et al.~\cite{jakob2019optimal} proposed a scalable multi-echelon $(s,Q)$ framework, equivalently referred to as the $(r,q)$ model with reorder point $r$ and order quantity $q$, that treats parking orbits as shared depots and selects policy parameters with an analytical cost model under Poisson failures and exponential launch lead times.}

\textcolor{black}{A complementary line of work considers dual-channel replenishment for constellation spares, motivated by terrestrial dual-sourcing studies~\cite{ramasesh1991sole,chiang1994sole,pan1991multiple,schimpel2014dual}. Related system-level studies have also shown that multiple deployment paths can provide structural benefits in constellation design~\cite{SUNG_deploy}.} Building on the multi-echelon $(s,Q)$ setting of Jakob et al.~\cite{jakob2019optimal}, Kim et al.~\cite{Kim2025_dual} extend the approach to two supply paths, \textcolor{black}{adopting the $(R_1,R_2,Q_1,Q_2)$ time-window policy from~\cite{schimpel2014dual}.} However, as formulated, this policy abstracts constellation-specific effects such as periodic alignment windows, \textcolor{black}{and therefore does not explicitly represent the inventory dynamics imposed by orbital mechanics.} In addition, the indirect path is fixed as the primary source, restricting the direct path to an auxiliary role\textcolor{black}{, which prevents the architecture from being channel-neutral}.

\textcolor{black}{This paper builds on several of our prior studies. Specifically, }our prior work introduced a Markov chain framework for constellation spare management and demonstrated an analytical evaluation on realistic scenarios~\cite{han2024analysis}. We then developed \textcolor{black}{channel-specific Markov-chain analyses for direct~\cite{han2025direct} and indirect~\cite{han2025indirect} replenishment, which capture small-LV lead-time stochasticity and parking-orbit drift dynamics, respectively}. This paper builds on that line \textcolor{black}{by presenting a Markov-chain-based framework for the early-phase analysis and design of hybrid spare-management architectures.} \textcolor{black}{The hybrid policy is modeled with a Markov chain that explicitly captures the periodic indirect replenishment opportunities and the stochastic direct deliveries, and computes the in-plane demand distribution from the model itself.} The model is channel-neutral\textcolor{black}{: it does not impose a preference between the direct and indirect paths,} and the optimization over policy parameters selects the preferred path and the primary and secondary roles\textcolor{black}{; with specific settings, the policy reduces to a purely direct or purely indirect case}. The method returns the full periodic stationary distribution for detailed cost and resilience metrics, and we \textcolor{black}{demonstrate its practical value in a real-world mega-constellation case study}.

The remainder of this paper is organized as follows. Section~\ref{sec2} introduces modeling preliminaries. Section~\ref{sec3} presents the cycle-based Markov analysis and an approximate analysis for the hybrid spare policy. Section~\ref{sec4} defines the cost model and performance measures. Section~\ref{sec5} validates the analysis with Monte Carlo simulations. Section~\ref{sec6} applies the framework to design optimization and a case study. Section~\ref{sec7} concludes the paper.

\section{Preliminaries} \label{sec2}
\subsection{Spare Management Policy}
\subsubsection{Hybrid Resupply Strategies}
The hybrid spare strategy combines both direct and indirect replenishment methods. In the direct method, a small launch vehicle is used to deliver spare satellites directly to an in-plane orbit when needed. \textcolor{black}{In the indirect method, a large launch vehicle sends spares to a parking orbit, where they wait for alignment in right ascension of the ascending node (RAAN) before being transferred to the target constellation orbit.} Previous studies have shown that the indirect method benefits from cost savings due to batch discounts, but the RAAN drift introduces longer delays compared to the direct method. On the other hand, the direct method provides faster response at a higher cost. \textcolor{black}{The hybrid strategy is considered as an architecture-level option that balances the lower launch cost of the indirect channel and the faster response of the direct channel.}

\begin{figure}[!h]
    \centering
    \includegraphics[width=.8\textwidth]{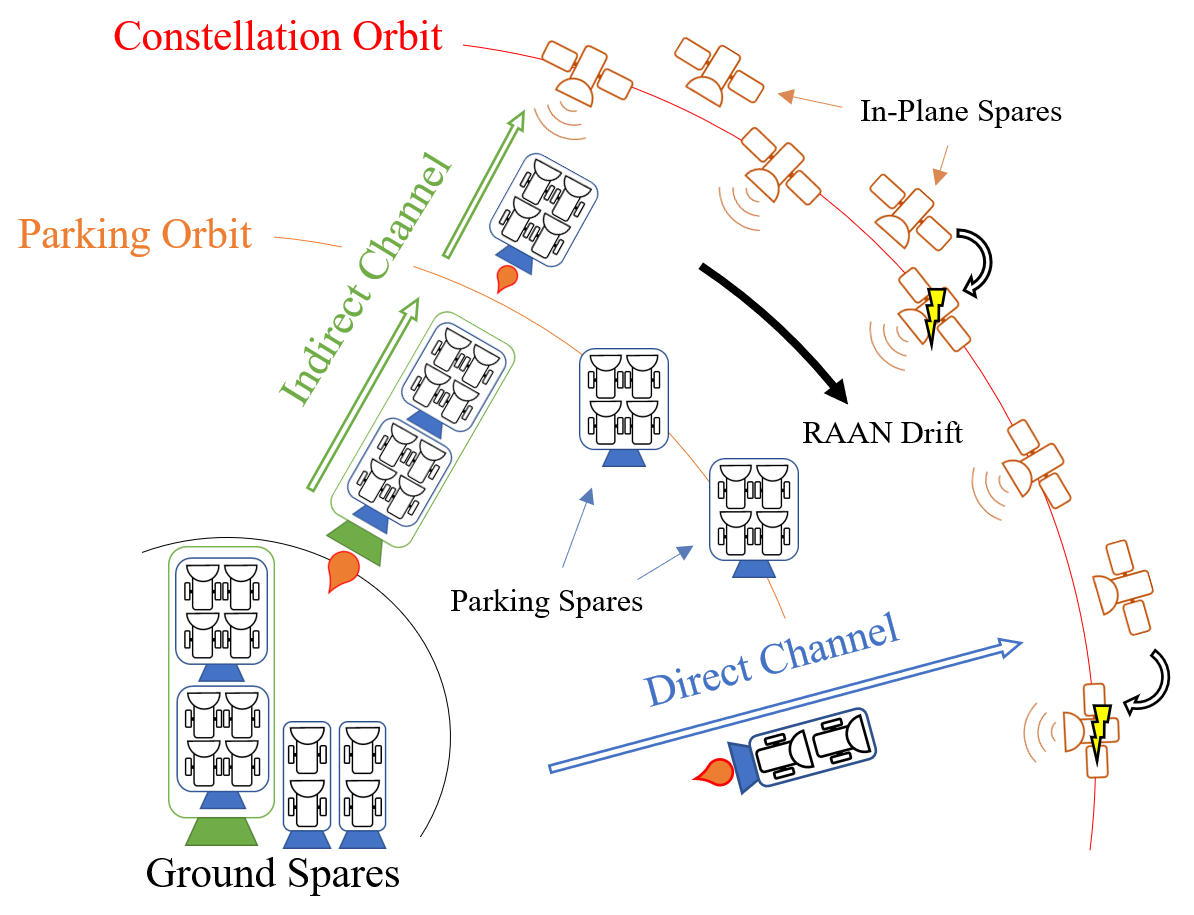}
    \caption{Illustration of Hybrid Spare Strategy}
    \label{fig:strategy}
\end{figure}

Figure~\ref{fig:strategy} illustrates the hybrid strategy. When a satellite fails, an in-plane spare is used for immediate replacement. \textcolor{black}{If the number of in-plane spares falls below a specified threshold, the policy can use either the indirect channel, by transferring spares from a parking orbit, or the direct channel, by launching new spares directly to the in-plane orbit.} At the same time, if the parking orbit stock drops below the threshold defined for parking orbits, a new resupply order is placed from the ground, and the replacement arrives after the launch vehicle lead time. In the indirect path, spare satellites are grouped into transfer buses and launched in batches using the large launch vehicle, allowing distribution to multiple constellation planes as needed.

\subsubsection{Inventory Management Policy} \label{sec:inventory_model}
In the hybrid strategy, constellation orbits follow a mixed dual-source reorder-point/order-quantity policy, with each channel governed by a different replenishment rule. The direct channel uses a standard $(r,q)$ policy, \textcolor{black}{in which an order of size $q$ is placed immediately when the stock level falls to or below the reorder point $r$, with no constraint on when orders may be placed.} The lead time is determined by the launch delay of the small launch vehicle, \textcolor{black}{which models the on-demand response of dedicated small-launch missions to in-plane orbits.}

\begin{figure}[!h]
    \centering
    \includegraphics[width=.75\textwidth, trim=5 5 5 5, clip]{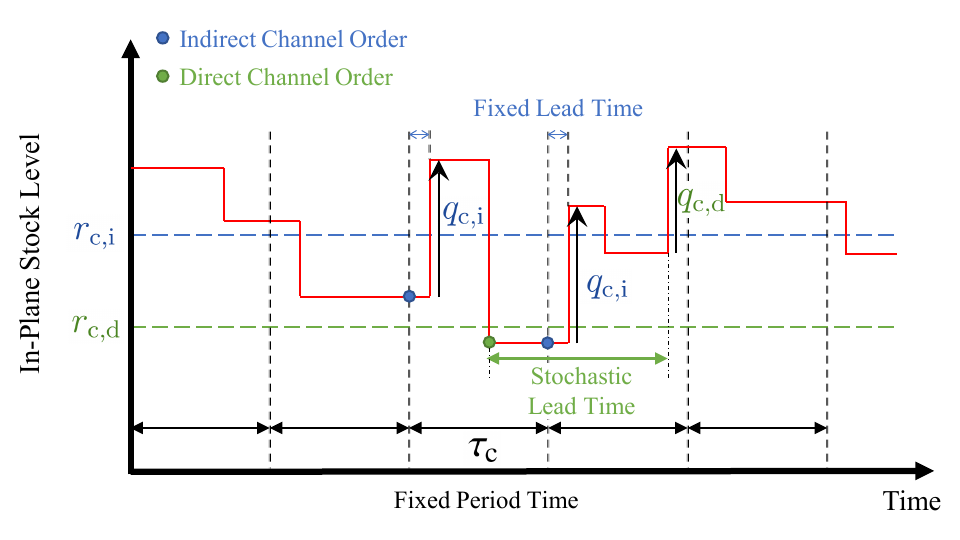}
    \includegraphics[width=.75\textwidth, trim=5 5 5 5, clip]{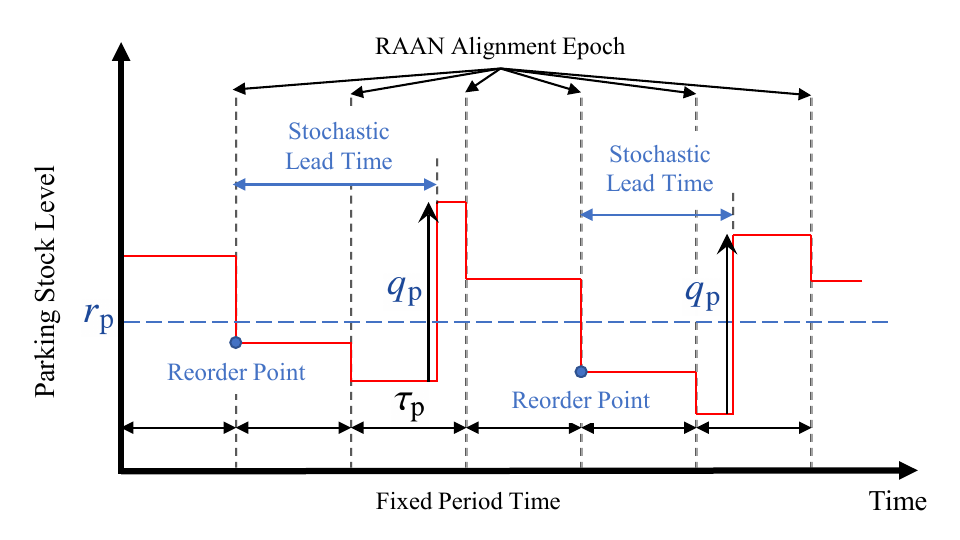}
    \caption{Stock level profile of (a) constellation and (b) parking orbits under hybrid strategy}
    \label{fig:inventory_model}
\end{figure}

\textcolor{black}{The indirect channel involves two stages, each governed by a periodic-review $(r,q,\tau)$ policy in which orders can only be placed at review epochs spaced by $\tau$. In our setting, $\tau$ is not a stock-monitoring interval as in classical inventory theory, but the interval between physically available replenishment opportunities imposed by orbital mechanics. The first stage replenishes the in-plane stock from parking orbits, with review epochs at successive RAAN alignments between a parking orbit and the constellation plane and a lead time equal to the fixed orbital-transfer duration. The second stage replenishes the parking stock from the ground, with review epochs at successive RAAN alignments between the parking orbit and constellation planes and a lead time determined by the heavy launch vehicle. See Fig.~\ref{fig:inventory_model} for more details.}

\subsection{Orbital Mechanics}
\subsubsection{Constellation Model}
We consider a large LEO Walker Delta constellation~\cite{walker1984satellite}. The system has $N_{\text{orbit}_\text{c}}$ circular planes with a common inclination and evenly spaced RAAN. Satellites in each plane are evenly phased in mean anomaly. Each plane carries $\bar N_{\text{sat}}$ nominal operational satellites to achieve the designed performance.

\subsubsection{Parking Orbits}
The indirect channel stages spares in $N_{\text{orbit}_\text{p}}$ circular parking orbits. These orbits use the same inclination as the constellation but a different semi-major axis, and they are evenly spaced in RAAN. The semi-major axis offset creates a differential RAAN drift that produces periodic alignment windows with the constellation planes.

\subsubsection{Orbital Maneuver}
Indirect transfers from parking to constellation planes use coplanar Hohmann maneuvers executed only at plane alignment. We treat a group of spares and a transfer bus as one stack during the maneuver. The fuel mass is
\begin{equation} \label{eq:fuel_mass}
    m_\text{fuel} = m_\text{dry}\!\left(\exp\!\left(\frac{\Delta v}{v_\text{ex}}\right) - 1\right),
\end{equation}
where $m_\text{dry}$ is the dry mass of the stack, $v_\text{ex}$ is the effective exhaust velocity, and $\Delta v$ is the Hohmann requirement \cite{prussing1993orbital}:
\begin{equation} \label{eq:delta_v}
    \Delta v =
    \sqrt{\frac{\mu_\oplus}{a_\text{p}}}\!\left(\sqrt{\frac{2a_\text{c}}{a_\text{p}+a_\text{c}}}-1\right)
    + \sqrt{\frac{\mu_\oplus}{a_\text{c}}}\!\left(1-\sqrt{\frac{2a_\text{p}}{a_\text{p}+a_\text{c}}}\right),
\end{equation}
with $\mu_\oplus$ the Earth gravitational parameter and $a_\text{p}$, $a_\text{c}$ the semi-major axes of the parking and constellation orbits ($a_\text{c}>a_\text{p}$). The Hohmann time of flight in LEO is on the order of hours and is neglected.

\subsubsection{RAAN Drift and Alignment}
Secular $J_2$ perturbations cause RAAN drift for circular orbits \cite{prussing1993orbital}:
\begin{equation} \label{eq:rann_drift}
    \dot \Omega_\text{c} = -\frac{3}{2} J_2 \frac{R_\oplus^2}{a_\text{c}^2}\sqrt{\frac{\mu_\oplus}{a_\text{c}^3}}\cos i,
    \qquad
    \dot \Omega_\text{p} = -\frac{3}{2} J_2 \frac{R_\oplus^2}{a_\text{p}^2}\sqrt{\frac{\mu_\oplus}{a_\text{p}^3}}\cos i,
\end{equation}
where $i$ is the common inclination, and $R_\oplus$, $\mu_\oplus$, $J_2$ are Earth parameters. Symmetric spacing implies periodic alignment between parking and constellation planes with intervals
\begin{equation} \label{eq:t_plane_indirect}
    \tau_\text{c} = \frac{2\pi}{N_{\text{orbit}_\text{p}}\lvert \dot\Omega_\text{c}-\dot\Omega_\text{p}\rvert},
    \qquad
    \tau_\text{p} = \frac{2\pi}{N_{\text{orbit}_\text{c}}\lvert \dot\Omega_\text{c}-\dot\Omega_\text{p}\rvert}.
\end{equation}
\textcolor{black}{These intervals define the RAAN-alignment replenishment opportunities used by the indirect channel.}

Based on the discrete Markov modeling in this work, we discretize time using a base interval $\tau_\text{mc}$. With a proper selection of $\tau_\text{mc}$ or through rounding operations, $\tau_\text{c}$ and $\tau_\text{p}$ can be expressed as integer multiples of $\tau_\text{mc}$ as
\begin{equation} 
    \tau_\text{p} = k_\text{p} \tau_\text{mc}, \qquad
    \tau_\text{c} = k_\text{c} \tau_\text{mc}, \qquad
    k_\text{p}, k_\text{c} \in \mathbb{N}.
\end{equation}
\textcolor{black}{The integers $k_\text{p}$ and $k_\text{c}$ define the RAAN-alignment cycle used in the policy analysis.}

\subsection{Markov Chain Model}
The spare satellite inventories in both constellation and parking orbits are modeled as discrete-time Markov chains. \textcolor{black}{Denoting the number of satellites, including spares, at time step $k$ by $X_k \in \{0,1,\ldots,N_\text{sat}\}$,} the state distribution is written as
\begin{equation}
    \pi_k =
    \begin{bmatrix}
        \mathbb{P}(X_k=N_\text{sat}) &
        \mathbb{P}(X_k=N_\text{sat}-1) &
        \cdots &
        \mathbb{P}(X_k=0)
    \end{bmatrix}^{\top},
\end{equation}
where the entries are ordered from state $N_\text{sat}$ down to state $0$ and satisfy $\mathbbm{1}\cdot\pi_k=1$.

For later use in the time-resolved formulation, we also define an unnormalized state distribution $\tilde{\pi}_k$. Its total weight is
\begin{equation}
    \gamma_k = \mathbbm{1}\cdot\tilde{\pi}_k \leq 1,
\end{equation}
and the corresponding normalized distribution is
\begin{equation}
    \pi_k = \frac{\tilde{\pi}_k}{\mathbbm{1}\cdot\tilde{\pi}_k}.
\end{equation}

Assuming a time-homogeneous transition matrix $P\in\mathbb{R}^{(N_\text{sat}+1)\times(N_\text{sat}+1)}$ with entries $P_{ij}=\mathbb{P}(X_{k+1}=i\mid X_k=j)$, the chain evolves as
\begin{equation}
    \pi_{k+1}=P\pi_k.
\end{equation}
Under standard ergodicity conditions, the Markov chain has a unique stationary distribution satisfying
\begin{equation}
    \pi=P\pi,
\end{equation}
which represents the long-run fraction of time spent in each state~\cite{hillier2015introduction}. In this spare-management model, failures decrease the state and replenishment increases the state.

\subsection{Probabilistic Model and Transition Matrix}
The probabilistic modeling framework used in this section, including the satellite failure distribution, lead time modeling, and transition matrix construction, is largely adopted from our previous work \cite{han2025indirect}.

\subsubsection{Satellite Failure Probability Distribution}
\textcolor{black}{This work models satellite failures as a Poisson process}, which is reasonable under the infinite-horizon stationary analysis adopted here, since fast in-orbit replacement yields an effective constant failure-rate model. \textcolor{black}{Denoting the failure rate of an operational satellite per $\tau_\text{mc}$ by $\lambda_\text{sat}$,} the probability of observing $k$ failures from $n$ satellites (including spares) during $\tau_\text{mc}$ is given by:
\begin{equation}\label{mpois_fail}
\begin{aligned}
    \nu_{k,n} &= 
    \mathbb{P}(F=k| X=n)\\
    &= \begin{cases}
        0 & \mbox{if } k > \bar{N}_\text{sat}\\ 
        \frac{(n\lambda_\text{sat})^k}{k!}e^{-n\lambda_\text{sat}} & \mbox{if } n \leq \bar{N}_\text{sat} \mbox{, } k \leq \bar{N}_\text{sat} \\ 
        \frac{(\bar{N}_\text{sat}\lambda_\text{sat})^k}{k!}e^{-\bar{N}_\text{sat}\lambda_\text{sat}} & \mbox{if } n > \bar{N}_\text{sat} \mbox{, } k \leq \bar{N}_\text{sat} \\
    \end{cases}
\end{aligned}
\end{equation}
where $F$ is the number of failures, $\bar{N}_\text{sat}$ is the nominal number of operational (non-spare) satellites. \textcolor{black}{An implicit assumption here is that failure replacement happens immediately and that spare satellites do not fail.} 
Furthermore, \textcolor{black}{denoting the maximum number of satellites including spares in constellation orbits by $N_{\text{sat}_\text{c}}$, the state transition matrix due to failure is given by}:
\begin{equation} \label{eq:fail_matrix_sim}
    P_{f} = \begin{bmatrix}
    \nu_{0,N_{\text{sat}_\text{c}}} 
      & 0                     & \cdots         & 0         \\
    \nu_{1,N_{\text{sat}_\text{c}}} 
      & \nu_{0,N_{\text{sat}_\text{c}}-1}   & \cdots         & 0         \\
    \vdots 
      & \vdots                & \ddots & \vdots       \\
    1-\sum_{k=0}^{N_{\text{sat}_\text{c}}-1} \nu_{k,N_{\text{sat}_\text{c}}}
      & 1-\sum_{k=0}^{N_{\text{sat}_\text{c}}-2} \nu_{k,N_{\text{sat}_\text{c}}-1} 
                             & \cdots 
                                       & \nu_{0,0}  
\end{bmatrix}
\end{equation}
where $P_f \in \mathbb{R}^{( N_{\text{sat}_\text{c}} +1 )\times ( N_{\text{sat}_\text{c}}+1)}$. By construction, each column vector sums to one, and the matrix is lower triangular, clearly showing that multiplying by $P_f$ always decreases the state level. \textcolor{black}{From an in-plane state distribution $\pi$, $P_f\pi$ gives the distribution after a one-step failure.}

\subsubsection{LV Lead-Time Probability Distribution}
\textcolor{black}{Similarly to our prior work~\cite{jakob2019optimal}, }the lead time for ground-based resupply is modeled using a shifted exponential distribution: 
\begin{equation*}
    T_{(\cdot)} \sim \text{Exp}(\mu_{\text{lv},(\cdot)}) + \tau_{\text{lv},(\cdot)}
\end{equation*}
where $(\cdot)$ denotes the channel type, either direct (d) or indirect (i). For simplicity, we drop the $(\cdot)$ notation in the following equations; the parameters $\mu_\text{lv}$ and $\tau_\text{lv}$ should be interpreted as referring to the specified channel. The corresponding probability density function is given by:
\begin{equation}
    f_\text{lv}(t; \mu_\text{lv}, \tau_\text{lv}) = 
    \begin{cases}
        \frac{1}{\mu_\text{lv}} e^{-(t - \tau_\text{lv}) / \mu_\text{lv}}, & t \geq \tau_\text{lv} \\
        0, & t < \tau_\text{lv}
    \end{cases}
\end{equation}
To simplify discrete-time modeling, the Markov chain time step $\tau_\text{mc}$ is chosen such that $\tau_\text{lv}$ is an integer multiple of $\tau_\text{mc}$, and we define
\begin{equation}
    m_\text{d} = \frac{\tau_\text{lv,d}}{\tau_\text{mc}},\quad 
    m_\text{i} = \frac{\tau_\text{lv,i}}{\tau_\text{mc}}.
\end{equation}
After the fixed delay $\tau_\text{lv}$ has elapsed, the residual lead time is exponentially distributed. By the memoryless property, the probability of receiving the order within the next $\tau_\text{mc}$, given that it has not yet arrived, is constant:
\begin{equation} \label{eq:rho}
    \rho = \mathbb{P}(T \leq \tau + \tau_\text{mc} \mid T > \tau) = 1 - e^{-\tau_\text{mc}/\mu_\text{lv}}, \quad \text{for any } \tau \geq \tau_\text{lv}.
\end{equation}
This allows the post-delay arrival process to be modeled as a geometric distribution with constant success probability $\rho$, assuming that each orbit places at most one LV order at a time. The direct and indirect channels accordingly have arrival probabilities $\rho_\text{d}$ and $\rho_\text{i}$, corresponding to $\mu_\text{lv,d}$ and $\mu_\text{lv,i}$.

\subsubsection{State Space and Reorder Threshold Projections}
To apply the \((r, q)\) policy, we need to isolate the portion of \(\pi\) corresponding to states where the stock level is less than or equal to \(r\). Let the maximum number of satellites (either in constellation or parking orbits) be \(N_{\text{sat}}\), so the state distribution \(\pi^{(\cdot)}\) lies in \(\mathbb{R}^{N_{\text{sat}} + 1}\). We define the following projection matrices:
\begin{equation} \label{eq:Cr_matrix}
\begin{aligned}
    C_{r}^+ = \begin{bmatrix}
        I_{N_{\text{sat}}-r} & \textbf{0}_{(N_{\text{sat}}-r)\times (r+1)} \\ \textbf{0}_{(r+1) \times (N_{\text{sat}}-r)} & \textbf{0}_{r+1}
    \end{bmatrix},\quad
    C_{r}^- = I_{N_\text{sat}+1} - C_{r}^+
\end{aligned}
\end{equation}
Then, \(C_{r}^+\pi\) gives the distribution for \(X > r\), and \(C_{r}^-\pi\) gives the distribution for \(X \leq r\). For example, the projection matrices for the in-plane and parking orbits are given by \(C_{r_\text{c,d}}^+, C_{r_\text{c,i}}^+ \in \mathbb{R}^{(N_{\text{sat}_\text{c}}+1)\times (N_{\text{sat}_\text{c}}+1)}\) and \(C_{r_\text{p}}^+ \in \mathbb{R}^{(N_{\text{sat}_\text{p}}+1)\times (N_{\text{sat}_\text{p}}+1)}\), respectively.

\subsubsection{Replenishment Transition Matrix}
After the lead time elapses, the system receives \(q_{(\cdot)}\) spare satellites. To model the corresponding state update, we define the replenishment transition matrix, which maps the state distribution immediately before delivery to the distribution immediately after. \textcolor{black}{For both the in-plane orbits served by the direct channel and the parking orbits served by the indirect channel}, the matrices $P_{q_\text{c,d}}$ and $P_{q_\text{p}}$ share the same structure and are defined as:
\begin{equation} \label{eq:resupply_trans_matrix}
    P_{q} =
        \begin{bmatrix}
        \begin{array}{c|c}
        I_q & I_{r+1} \\
        \mathbf{0}_{(r+1)\times q} & \mathbf{0}_{q\times (r+1)}
        \end{array}
        \end{bmatrix},
\end{equation}
where \(P_{q} \in \mathbb{R}^{(N_{\text{sat}} + 1) \times (N_{\text{sat}} + 1)}\) under $(r,q)$ policy. \textcolor{black}{The block \(I_q\) leaves states above the reorder point unchanged, while \(I_{r+1}\) maps states with \(X\le r\) to \(X+q\). This keeps \(P_q\) column-stochastic over the full state space.}

Specifically, \(P_{q_\text{c,d}}\) and \(P_{q_\text{p}}\) are the replenishment matrices for the in-plane (direct) and parking orbits, respectively, with dimensions \((N_{\text{sat}_\text{c}}+1) \times (N_{\text{sat}_\text{c}}+1)\) and \((N_{\text{sat}_\text{p}}+1) \times (N_{\text{sat}_\text{p}}+1)\). In contrast, the matrix \(P_{q_\text{c,i}}\), which models indirect in-plane replenishment, has a different structure and is described separately in Eq.~\eqref{eq:park_avail_mat}.

\section{Modeling and Analysis of Spare Management Policy} \label{sec3}
The analysis approach for the hybrid resupply strategy follows the general structure of our previous work on the indirect method \cite{han2025indirect}. In this approach, the constellation-orbit and parking-orbit Markov chains are modeled separately, with consistency enforced between them. Although the overall framework is similar to the previous work, the detailed formulation differs substantially to accommodate the hybrid strategy. We present each component of the model step by step in the following subsections.

Throughout this section, subscripts $(\cdot)_\text{c}$ and $(\cdot)_\text{p}$ refer to the constellation and parking orbits, respectively. When it is necessary to distinguish between the direct and indirect channels of the constellation orbits, we use the subscripts $(\cdot)_\text{c,d}$ and $(\cdot)_\text{c,i}$, respectively (e.g., $q_\text{c,d}$ and $q_\text{c,i}$).

\subsection{Cycle-Based Markov Modeling Framework} \label{sec:event_vs_cycle}
In our previous work on the indirect strategy \cite{han2025indirect}, we adopted an event-based approach, which models inventory transitions between key events, such as when an order is placed ($\pi^q$) and when it arrives ($\pi^r$). This formulation yields compact analytical expressions and requires fewer state variables. However, it also necessitates implicitly accounting for all possible transitions during the intervals between events, which makes the derivation nontrivial, particularly in a hybrid strategy involving dual replenishment channels.

In contrast, the present work adopts a time-resolved approach, where the inventory distribution is explicitly modeled at each discrete time step. The state is tracked during the inter-order period ($\pi^{\text{io}}$) and throughout the lead-time phase ($\pi^{\text{lt}}_m$). While this increases the state space dimensionality, it simplifies transition modeling and is well suited for hybrid replenishment structures. The assumption of a memoryless lead time distribution further ensures a finite state space.

\begin{figure}[!h]
    \centering
    \includegraphics[width=.8\textwidth]{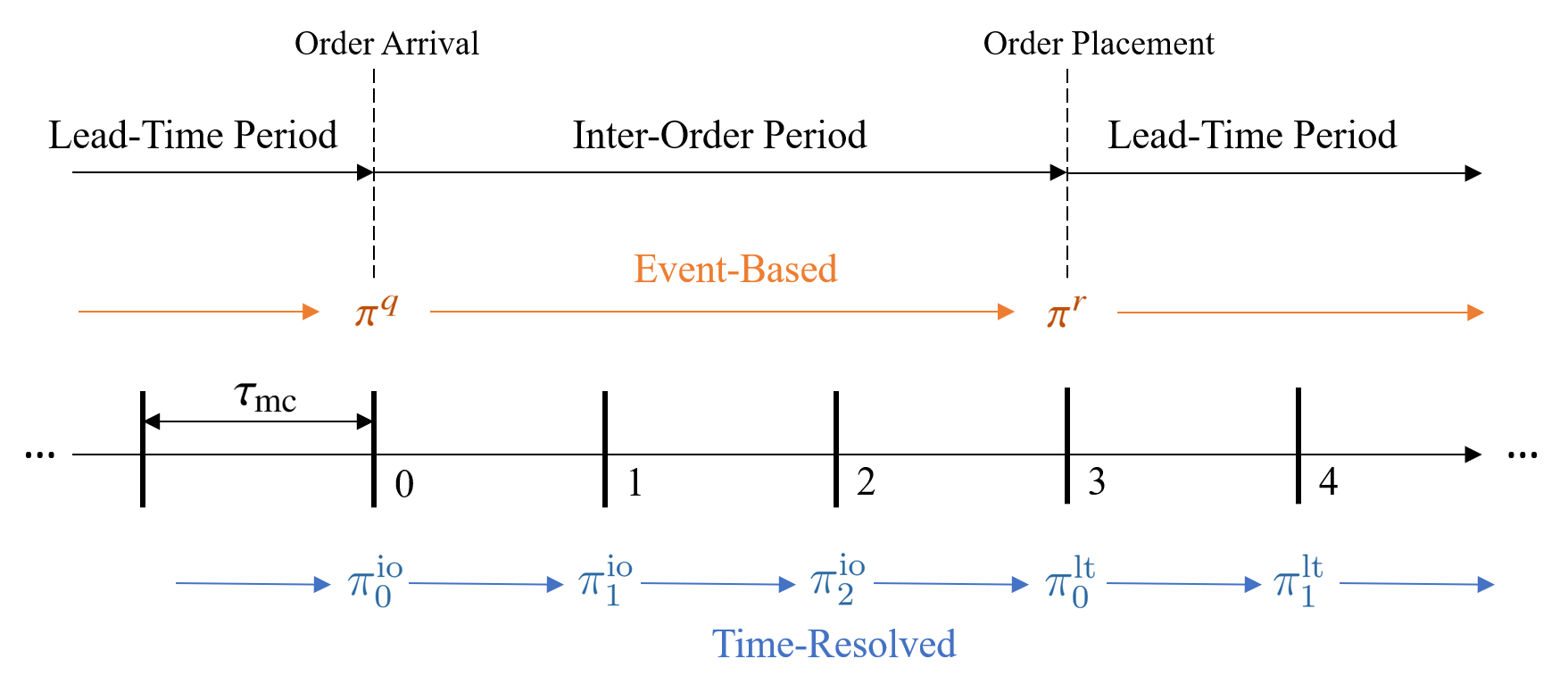}
    \caption{Illustration of Two Approaches}
    \label{fig:method_diff}
\end{figure}

Although both approaches are mathematically equivalent when transitions are properly accounted for, they differ in the trade-off between model compactness and derivation complexity. Given the structure of the hybrid system considered here and the use of shifted exponential lead times, the time-resolved approach offers greater flexibility and ease of implementation, at the cost of added dimensionality. The conceptual differences between the two methods are illustrated in Fig.~\ref{fig:method_diff}.

Finally, the inter-order (IO) and lead-time (LT) phases are defined independently for each layer based on their respective ordering conditions. For the constellation orbits, the phases are determined by the direct channel using small launch vehicles, while the indirect channel acts as a fixed periodic input. Conversely, for the parking orbits, the IO and LT phases are defined by the indirect channel using large launch vehicles.

\subsection{Constellation Orbit Analysis Method} \label{sec:in_plane_anal_indirect}
In this section, we introduce the method for computing the periodic stationary distribution of in-plane inventory.
First, the maximum number of satellites in a constellation orbit is
\begin{equation*}
    N_{\text{sat}_\text{c}}
    = \max\left\{r_{\text{c,d}} + q_{\text{c,d}},\;
    r_{\text{c,i}} + q_{\text{c,i}} + q_{\text{c,d}}\right\},
\end{equation*}
in units of satellites. This bound captures the worst-case scenario in which a direct batch arrives immediately after an indirect replenishment.

\subsubsection{Indirect Replenishment Transition Matrix} \label{sec:in_plane_replenishment_indirct}
We assume that the parking-orbit inventory \(X_\text{p}\) is identically and independently distributed (i.i.d), and that its distribution just before a RAAN contact, \(\mathbb{P}(X_\text{p}=k \mid E)\), is known (see Section~\ref{sec:park_avail}). The parking availability probability is defined as
\begin{equation} \label{eq:parking_avail}
    \kappa_j = \mathbb{P}(X_\text{p} \geq j \mid E),
\end{equation}
where \(E\) denotes the RAAN contact event, and \(j = D_\text{c}\) is the realized in-plane demand expressed in integer batches of size \(q_\text{c,i}\), given by
\begin{equation} \label{eq:park_demand_func}
    D_\text{c} = \begin{cases}
        \left\lceil \dfrac{r_\text{c,i} + 1 - X_\text{c}}{q_\text{c,i}} \right\rceil, & \text{if } X_\text{c} \leq r_\text{c,i}, \\
        0, & \text{otherwise},
    \end{cases}
\end{equation}
where \(\lceil \cdot \rceil\) is the ceiling function. If \(X_\text{c} > r_\text{c,i}\), no replenishment is needed. Otherwise, \(D_\text{c}\) gives the number of batches required to raise the inventory above the reorder point. Thus, \(\kappa_j\) is the probability that at least \(j\) batches are available in the parking orbit at the RAAN alignment.

The transition matrix for indirect replenishment is then defined as
\begin{equation} \label{eq:park_avail_mat}
    P_{q_\text{c,i}} = 
\begin{bmatrix}
    \kappa_0 I_{w} & \kappa_1 I_{q_\text{c,i}} & \kappa_2 I_{q_\text{c,i}} & \cdots \\
    0 & (\kappa_0 - \kappa_1) I_{q_\text{c,i}} & (\kappa_1 - \kappa_2) I_{q_\text{c,i}} & \cdots \\
    0 & 0 & (\kappa_0 - \kappa_1) I_{q_\text{c,i}} & \cdots \\
    \vdots & \vdots & \vdots & \ddots
\end{bmatrix},
\end{equation}
where \(w=N_{\text{sat}_\text{c}}-r_\text{c,i}\), and \(P_{q_\text{c,i}} \in \mathbb{R}^{(N_{\text{sat}_\text{c}}+1)\times (N_{\text{sat}_\text{c}}+1)}\). Each block  corresponds to a demand level and the probability of fulfilling it.

For example, the first \(w\) entries of the state vector \(\pi\) (i.e., \(r_\text{c,i} < X_\text{c} \leq N_{\text{sat}_\text{c}}\)) require no replenishment and are multiplied by \(\kappa_0 I_w\). The next \(q_\text{c,i}\) entries (i.e., \(r_\text{c,i} - q_\text{c,i} < X_\text{c} \leq r_\text{c,i}\)) represent a one-batch demand. If a spare is available with probability \(\kappa_1\), the inventory increases by \(q_\text{c,i}\), which corresponds to the (1,2) block; otherwise, the state remains unchanged with probability \(\kappa_0 - \kappa_1\), as shown in the (2,2) block.

In summary, if \(\pi\) denotes the in-plane inventory distribution just before RAAN contact, then \(P_{q_\text{c,i}} \pi\) yields the distribution after accounting for parking availability.

\subsubsection{Periodic Transition Matrix of Constellation Orbit } \label{sec:transition_const}
\begin{figure}[!h]
    \centering
    \includegraphics[width=.8\textwidth]{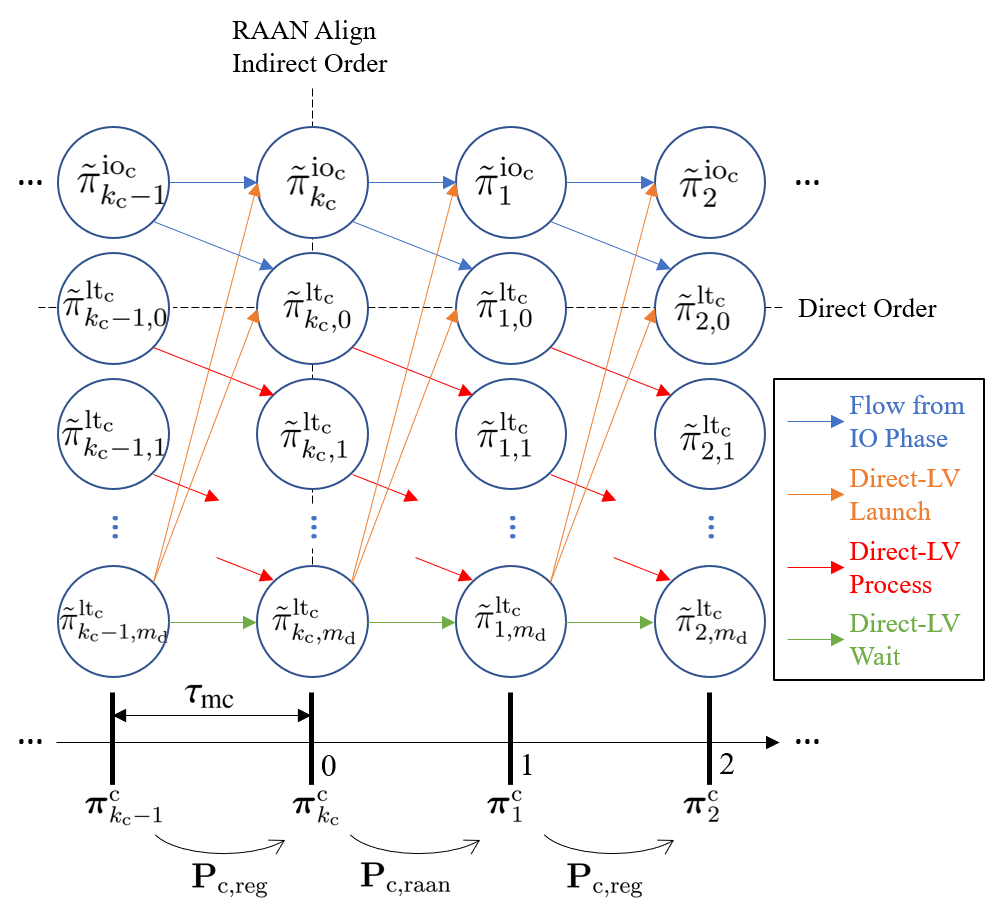}
    \caption{In-plane state transitions under hybrid replenishment, showing IO/LT phases and direct-LV flow paths.}
    \label{fig:transition_diagram}
\end{figure}

In this subsection, we derive the transition matrix between each pair of weighted distributions \(\tilde{\pi}^{\text{io}_\text{c}}_k\) and \(\tilde{\pi}^{\text{lt}_\text{c}}_{k,m}\). Refer to Fig.~\ref{fig:transition_diagram} for overall transition between state distributions.

\paragraph{Transition into IO Phase}
The transition into the IO phase at time step \(k+1\) may occur in two ways: the system either remains in the IO phase from time \(k\), or it exits the LT phase after completing the constant processing delay (\(m = m_\text{d}\)) and receives a direct replenishment (\(q_\text{c,d}\)) with probability \(\rho_\text{d}\), which arises from the memoryless property of the lead-time distribution. These two sources of inflow are combined into a single pre-update distribution:
\[
    \tilde{\pi}^{\text{io}_\text{c}}_k + \rho_\text{d} P_{q_\text{c,d}} \tilde{\pi}^{\text{lt}_\text{c}}_{k,m_\text{d}}.
\]
Failures may occur between steps \(k\) and \(k+1\), so the combined distribution is updated as:
\[
    P_f \left( \tilde{\pi}^{\text{io}_\text{c}}_k + \rho_\text{d} P_{q_\text{c,d}} \tilde{\pi}^{\text{lt}_\text{c}}_{k,m_\text{d}} \right).
\]

Using this pre-update distribution, the final IO-phase distribution at time \(k+1\) is given by
\begin{equation} \label{eq:to_io_cases}
\tilde{\pi}^{\text{io}_\text{c}}_{k+1} = 
\begin{cases}
C^+_{r_\text{c,d}} \left( P_f \tilde{\pi}^{\text{io}_\text{c}}_k + \rho_\text{d} P_f P_{q_\text{c,d}} \tilde{\pi}^{\text{lt}_\text{c}}_{k,m_\text{d}} \right),
& k = 1,\dots,k_\text{c}-1, \\
C^+_{r_\text{c,d}} \left[ 
P_{q_\text{c,i}} \left( P_f \tilde{\pi}^{\text{io}_\text{c}}_{k_\text{c}} + \rho_\text{d} P_f P_{q_\text{c,d}} \tilde{\pi}^{\text{lt}_\text{c}}_{k_\text{c},m_\text{d}} \right) 
\right], 
& k = k_\text{c},
\end{cases}
\end{equation}
In the first case, no RAAN alignment occurs, so only the direct replenishment condition is checked. In the second case, indirect replenishment is applied first via \(P_{q_\text{c,i}}\), followed by direct policy filtering via \(C^+_{r_\text{c,d}}\).

\paragraph{Transition into and within the LT Phase}
The transition into a new LT phase, i.e., \(m=0\), at time step \(k+1\) occurs in the same way as the transition into the IO phase. Specifically, the following pre-update distribution represents the portion of the system that newly enters the LT phase:
\[
     \tilde{\pi}^{\text{io}_\text{c}}_k + \rho_\text{d} P_{q_\text{c,d}} \tilde{\pi}^{\text{lt}_\text{c}}_{k,m_\text{d}}.
\]
The resulting distribution for \(m = 0\) is then given by:
\begin{equation}
\tilde{\pi}^{\text{lt}_\text{c}}_{k+1,\,0}
=
\begin{cases}
  C^-_{r_\text{c,d}} \left( P_f\,\tilde{\pi}^{\text{io}_\text{c}}_k + \rho_\text{d} P_f P_{q_\text{c,d}} \tilde{\pi}^{\text{lt}_\text{c}}_{k,m_\text{d}} \right),
    &  k = 1,\dots,k_\text{c}-1, \\
  C^-_{r_\text{c,d}} \left[ P_{q_\text{c,i}} \left( P_f\,\tilde{\pi}^{\text{io}_\text{c}}_k + \rho_\text{d} P_f P_{q_\text{c,d}} \tilde{\pi}^{\text{lt}_\text{c}}_{k,m_\text{d}} \right) \right],
    &  k = k_\text{c}.
\end{cases}
\end{equation}

Next, the transition from LT to LT for \(m < m_\text{d}\) is straightforward. The pre-update distribution in this case is:
\[
 \tilde{\pi}^{\text{lt}_\text{c}}_{k,m},
\]
and the update equation becomes:
\begin{equation}
\tilde{\pi}^{\text{lt}_\text{c}}_{k+1,\,m+1}
=
\begin{cases}
  P_f\,\tilde{\pi}^{\text{lt}_\text{c}}_{k,m}, 
    &   k = 1,\dots,k_\text{c}-1, \\
  P_{q_\text{c,i}}\,P_f\,\tilde{\pi}^{\text{lt}_\text{c}}_{k,m}, 
    &  k = k_\text{c},
\end{cases}
\end{equation}
for $m < m_\text{d}$.

Finally, when \(m = m_\text{d}\), the distribution has two sources: a portion newly enters from \(m_\text{d}-1\), and the rest remains at \(m_\text{d}\) due to the memoryless property of the lead time. The combined pre-update distribution is:
\[
     \tilde{\pi}^{\text{lt}_\text{c}}_{k,m_\text{d}-1} + (1-\rho_\text{d}) \tilde{\pi}^{\text{lt}_\text{c}}_{k,m_\text{d}},
\]
and the corresponding update equation is:
\begin{equation}
\tilde{\pi}^{\text{lt}_\text{c}}_{k+1,\,m_\text{d}}
=
\begin{cases}
  P_f\,\tilde{\pi}^{\text{lt}_\text{c}}_{k,\,m_\text{d}-1}
  + (1-\rho_\text{d})\,P_f\,\tilde{\pi}^{\text{lt}_\text{c}}_{k,\,m_\text{d}},
    & k = 1,\dots,k_\text{c}-1, \\
  P_{q_\text{c,i}} \left( P_f\,\tilde{\pi}^{\text{lt}_\text{c}}_{k,\,m_\text{d}-1}
  + (1-\rho_\text{d})\,P_f\,\tilde{\pi}^{\text{lt}_\text{c}}_{k,\,m_\text{d}} \right),
    & k = k_\text{c}.
\end{cases}
\end{equation}

\paragraph{Full Transition between IO and LT Phase}
Based on the update equations, we construct the full transition matrix by concatenating all sub-state distributions into a single vector:
\begin{equation}
    \boldsymbol{\pi}^\text{c}_k = \left(\tilde{\pi}^{\text{io}_\text{c}}_{k},\
    \tilde{\pi}^{\text{lt}_\text{c}}_{k,\,0},\,\dots,\,
    \tilde{\pi}^{\text{lt}_\text{c}}_{k,\,m_\text{d}}\right),
\end{equation}
where \( \mathbbm{1} \cdot \boldsymbol{\pi}^\text{c}_k = 1 \). The transition matrix for the case without RAAN alignment (\(k \neq k_\text{c}\), i.e., regular steps) is defined as
\begin{equation} \label{eq:P_c_reg_mat}
    \mathbf{P}_{\text{c},\text{reg}} = 
    \begin{bmatrix}
        C^+_{r_\text{c,d}}P_{f} & 0  & \cdots & 0 & \rho_\text{d} C^+_{r_\text{c,d}} P_{f} P_{q_\text{c,d}} \\
        C^-_{r_\text{c,d}}P_{f} & 0  & \cdots  & 0 & \rho_\text{d} C^-_{r_\text{c,d}} P_{f} P_{q_\text{c,d}} \\
        0 & P_{f} & \cdots & 0 & 0 \\
        \vdots & \vdots & \ddots  & \vdots & \vdots \\
        0 & 0 &  \cdots &   P_{f} & (1-\rho_\text{d})P_{f}
    \end{bmatrix}
\end{equation}
When RAAN alignment occurs (\(k = k_\text{c}\)), the transition matrix becomes:
\begin{equation} \label{eq:P_c_raan_mat}
    \mathbf{P}_{\text{c},\text{raan}} = 
    \begin{bmatrix}
        C^+_{r_\text{c,d}} P_{q_\text{c,i}}P_{f} & 0  & \cdots  & 0 & \rho_\text{d} C^+_{r_\text{c,d}} P_{q_\text{c,i}} P_{f} P_{q_\text{c,d}} \\
        C^-_{r_\text{c,d}} P_{q_\text{c,i}}P_{f} & 0  & \cdots  & 0 & \rho_\text{d} C^-_{r_\text{c,d}} P_{q_\text{c,i}} P_{f} P_{q_\text{c,d}} \\
        0 & P_{q_\text{c,i}}P_{f} & \cdots  & 0 & 0 \\
        \vdots & \vdots & \ddots  & \vdots & \vdots \\
        0 & 0 &  \cdots &  P_{q_\text{c,i}}P_{f} & (1-\rho_\text{d})P_{q_\text{c,i}}P_{f}
    \end{bmatrix}
\end{equation}
Here, \( \mathbf{P}_{\text{c},\text{reg}}, \mathbf{P}_{\text{c},\text{raan}} \in \mathbb{R}^{(m_\text{d}+2)(N_{\text{sat}_\text{c}}+1) \times (m_\text{d}+2)(N_{\text{sat}_\text{c}}+1)} \).

Using this formulation, the update rule for \(\boldsymbol{\pi}^\text{c}_k\) becomes:
\begin{equation} \label{eq:update_eq_constel}
    \boldsymbol{\pi}^\text{c}_{k+1} = 
    \begin{cases}
        \mathbf{P}_{\text{c},\text{reg}}\,\boldsymbol{\pi}^{\text{c}}_{k}, 
        & \text{if } k = 1,\dots,k_\text{c}-1, \\
        \mathbf{P}_{\text{c},\text{raan}}\,\boldsymbol{\pi}^{\text{c}}_{k}, 
        & \text{if } k = k_\text{c}.
    \end{cases}
\end{equation}
Finally, the periodic steady-state condition at $k=1$ is expressed as:
\begin{equation} \label{eq:stationary_constel}
\begin{aligned}
    \boldsymbol{\pi}^\text{c}_{1} 
    &= \mathbf{P}_{\text{c},\text{raan}} \left(\mathbf{P}_{\text{c},\text{reg}} \right)^{k_\text{c}-1} \boldsymbol{\pi}^\text{c}_{1} \\
    &= \mathbf{P}_{\text{c},\text{trn}} \boldsymbol{\pi}^\text{c}_{1},
\end{aligned}
\end{equation}
which serves as a fixed-point equation to determine the long-run state distribution of the constellation orbit inventory over a full RAAN alignment cycle.

\begin{remark}
\(C^+_{r_\text{c,d}}(\cdot)\) and \(C^-_{r_\text{c,d}}(\cdot)\) represent the direct \((r_\text{c,d}, q_\text{c,d})\) policy, while \(P_{q_\text{c,i}}(\cdot)\) corresponds to the indirect \((r_\text{c,i}, q_\text{c,i})\) policy triggered by RAAN alignment. The order of matrix reflects the assumed sequence within a time step: direct replenishment from a previous order first \((P_{q_\text{c,d}})\), followed by failures \((P_{f})\), then indirect replenishment due to RAAN alignment \((P_{q_\text{c,i}})\), and finally direct reorder \((C^-_{r_\text{c,d}})\), i.e., \(C^-_{r_\text{c,d}} P_{q_\text{c,i}} P_{f} P_{q_\text{c,d}}\). If this sequence changes, the corresponding matrix operations must be adjusted accordingly.
\end{remark}

\begin{remark}
The index \(k\) is treated cyclically with period \(k_\text{c}\), i.e.,
\(\tilde{\pi}^{\text{io}_\text{c}}_{k_\text{c}+i} = \tilde{\pi}^{\text{io}_\text{c}}_i,\
\tilde{\pi}^{\text{lt}_\text{c}}_{k_\text{c}+i,m} = \tilde{\pi}^{\text{lt}_\text{c}}_{i,m},\ \forall i \in \mathbb{Z}_{\ge 0}\).
This reflects the periodic RAAN alignment between constellation and parking orbits. As the orbital configuration repeats every \(k_\text{c}\) steps, the Markov model evolves over a finite cycle, allowing steady-state analysis within one period.
\end{remark}

\subsubsection{Average Distribution over the Replenishment Cycle}
By solving the linear system in Eq.~\eqref{eq:stationary_constel}, we obtain the stationary distribution $\boldsymbol{\pi}^\text{c}_{k}$ at each step $k$ of the review period. The IO and LT phases are mutually exclusive and together span the entire replenishment cycle. The weight $\gamma_k$ associated with each $\tilde{\pi}_k$ represents the fraction of time the system spends in that specific phase at step $k$, similar to how a stationary distribution in a standard Markov chain reflects the long-run frequency of visiting each state. 

Let \( \mathbf{I}_\text{c} \) denote the row block matrix that sums across all sub-distributions in the full state vector \( \boldsymbol{\pi}^{\text{c}}_k \), defined as
\begin{equation}
    \mathbf{I}_\text{c} = \begin{bmatrix} I_{N_{\text{sat}_\text{c}} + 1} & I_{N_{\text{sat}_\text{c}} + 1} & \cdots & I_{N_{\text{sat}_\text{c}} + 1} \end{bmatrix},
\end{equation}
where the total number of identity blocks is \( m_\text{d} + 2 \), corresponding to the IO phase and all \( m_\text{d} + 1 \) LT sub-phases. The operator \( \mathbf{I}_\text{c} \boldsymbol{\pi} \) marginalizes the full state vector to obtain the distribution of the constellation stock level.

With this definition, the average distribution at step \( k \) of the review period is given by
\begin{equation}
    \pi^{\text{rc}_\text{c}}_k = \mathbf{I}_\text{c} \boldsymbol{\pi}^{\text{c}}_k = 
    \tilde{\pi}^{\text{io}_\text{c}}_{k} + \sum_{m=0}^{m_\text{d}} \tilde{\pi}^{\text{lt}_\text{c}}_{k,m},
\end{equation}
and the average distribution over the entire replenishment cycle becomes
\begin{equation}
    \pi^{\text{rc}_\text{c}} = \frac{1}{k_\text{c}} \sum_{k=1}^{k_\text{c}} \pi^{\text{rc}_\text{c}}_k.
\end{equation}

\subsubsection{Computation of Average Replenishment Period}
To compute how frequently the direct channel is used, we calculate the average time interval for one replenishment cycle. First, define the following constants:
\begin{equation}
    c^{\text{io}_\text{c}} = \sum_{k=1}^{k_\text{c}} \sum_{i=0}^{N_{\text{sat}_\text{c}}} \tilde{\pi}^{\text{io}_\text{c}}_k(i)
    = \sum_{k=1}^{k_\text{c}} \gamma^{\text{io}_\text{c}}_k,
\end{equation}
\begin{equation}
    c^{\text{lt}_\text{c}} = \sum_{k=1}^{k_\text{c}} \sum_{m=0}^{m_\text{d}} \sum_{i=0}^{N_{\text{sat}_\text{c}}} \tilde{\pi}^{\text{lt}_\text{c}}_{k,m}(i) 
    = \sum_{k=1}^{k_\text{c}} \sum_{m=0}^{m_\text{d}} \gamma^{\text{lt}_\text{c}}_{k,m}.
\end{equation}
From the lead-time model, the average duration of the LT phase is given by
\begin{equation}
    \tau_{\text{lt}_\text{c}} = m_\text{d} \tau_\text{mc} + \mu_\text{lv,d}.
\end{equation}
Using this and the previously defined constants, the average duration of the IO phase is computed as
\begin{equation}
    \tau_{\text{io}_\text{c}} = \frac{c^{\text{io}_\text{c}}}{c^{\text{lt}_\text{c}}} \tau_{\text{lt}_\text{c}}.
\end{equation}
Therefore, the average time interval for one replenishment cycle via the direct channel is given by
\begin{equation}
    \tau_{\text{rc}_\text{c}} = \tau_{\text{lt}_\text{c}} + \tau_{\text{io}_\text{c}}.
\end{equation}

\subsubsection{Demand Distribution of Constellation Orbits to Parking Orbits} \label{sec:dmd_dist}
This section derives the probability distribution of the demand placed on the parking orbit. This demand is generated by the constellation orbits during each review period \(\tau_\text{c}\) and is expressed in units of the batch size \(q_\text{c,i}\).

The constellation-orbit inventory distribution immediately before indirect replenishment at the RAAN-alignment event $E$ is denoted by $\pi^{\text{rc}_\text{c}|E}$, or equivalently by $\pi^{r_\text{c,i}}$ when emphasizing the indirect-channel pre-replenishment state:
\begin{equation} \label{eq:pi_r_c_i}
    \pi^{\text{rc}_\text{c}|E}
    \equiv
    \pi^{r_\text{c,i}}
    =
    \mathbf{I}_\text{c}\,
    \mathbf{P}_\text{c,reg}\,
    \boldsymbol{\pi}^{\text{c}}_{k_\text{c}}.
\end{equation}
\textcolor{black}{Here, $\mathbf{P}_\text{c,reg}$ applies the regular one-step effects before indirect replenishment, such as direct-channel arrivals and failures. Since $\mathbf{I}_\text{c}$ sums the IO/LT components, the $C^+_{r_\text{c,d}}/C^-_{r_\text{c,d}}$ branching inside $\mathbf{P}_\text{c,reg}$ does not change the extracted marginal stock distribution.}

Next, we map this satellite inventory level to a batch demand \(D_\text{c}\) using Eq.~\eqref{eq:park_demand_func}. Since the demand is quantized by the batch size, multiple inventory states may result in the same demand. For example, if the deficit is 1 or 2 satellites and the batch size is 2, both states trigger a demand of 1 batch. To compute the probability mass function \(\chi_j\) for demanding exactly \(j\) batches, define
\begin{equation}
    \mathcal{S}_j =
    \left\{x \in \{0,1,\dots,N_{\text{sat}_\text{c}}\}: D_\text{c}(x)=j \right\},
\end{equation}
where \(D_\text{c}(x)\) denotes Eq.~\eqref{eq:park_demand_func} evaluated at \(X_\text{c}=x\). Then,
\begin{equation} \label{eq:inplane_demand_prob}
\begin{aligned}
    \chi_j
    &= \mathbb{P}(D_\text{c}=j \mid E) \\
    &= \sum_{x\in\mathcal{S}_j}  \pi^{\text{rc}_\text{c}|E}(X_\text{c}=x),
    \quad j=0,1,\dots,\left\lceil\frac{r_\text{c,i}+1}{q_\text{c,i}}\right\rceil .
\end{aligned}
\end{equation}
The resulting \(\chi\)-vector serves as the demand input for the parking-orbit analysis.

For example, if $r_\text{c,i} = 3$, $q_\text{c,i} = 2$, and $N_{\text{sat}_\text{c}} = 5$, the states $X_\text{c} \in \{4,5\}$ are above the reorder point and therefore correspond to zero demand. Thus, $\chi_0$ is the sum of their probabilities:
\begin{equation*}
\chi_0 = \pi^{\text{rc}_\text{c}|E}(X_\text{c} = 4) + \pi^{\text{rc}_\text{c}|E}(X_\text{c} = 5).
\end{equation*}

Finally, the state distribution of the constellation orbit after the indirect replenishment opportunity, accounting for parking availability, is computed as
\begin{equation}
    \pi^{q_\text{c,i}}
    =
    P_{q_\text{c,i}}\,\pi^{r_\text{c,i}}
    =
    P_{q_\text{c,i}}\,\pi^{\text{rc}_\text{c}|E},
\end{equation}
where $P_{q_{\text{c,i}}}$ is the replenishment transition matrix defined in Section~\ref{sec:in_plane_replenishment_indirct}.

\subsection{Parking Orbit Analysis Method} \label{sec:event_park_analysis}
In this subsection, we present the method for computing the stationary distribution of the parking orbit inventory. The replenishment model for the parking inventory is identical to that in our previous work, and thus the method described in \cite{han2025indirect}, which adopts an event-based approach, can be directly applied. This remains the preferred method in our current model due to its efficiency and simplicity. However, for completeness of the proposed framework, we also provide a time-resolved formulation of the parking inventory model, which can be easily extended to incorporate additional complexities. Lastly, the maximum parking-orbit state level is given by
\begin{equation*}
N_{\text{sat}_\text{p}} = r_{\text{p}} + q_{\text{p}},
\end{equation*}
measured in units of batches.

\subsubsection{Failure and Replenishment Transition Matrix of Parking Orbits}
Assuming that spare satellites do not fail, the number of parking spares decreases only when they are transferred to the constellation orbits to meet demand. Using the demand distribution from Eq.~\eqref{eq:inplane_demand_prob}, the (demand-induced) failure transition matrix of the parking orbit at RAAN contact is
\begin{equation} \label{eq:fail_matrix_park}
    P_{f_\text{p}} =
    \begin{bmatrix}
        \chi_{0} & 0 & \cdots & 0\\
        \chi_{1} & \chi_{0} & \cdots & 0 \\
        \vdots & \vdots & \ddots & \vdots \\
        1-\sum_{i=0}^{N_{\text{sat}_\text{p}}-1}\chi_i & 1-\sum_{i=0}^{N_{\text{sat}_\text{p}}-2}\chi_i & \cdots & 1
    \end{bmatrix},
\end{equation}
which has the same lower triangular structure as $P_f$ in Eq.~\eqref{eq:fail_matrix_sim}, with the failure probabilities $\nu$ replaced by the (demand-induced) failure probabilities $\chi$. In summary, letting $\pi$ be the parking state distribution before RAAN alignment, the product $P_{f_\text{p}}\pi$ gives the distribution immediately after alignment (i.e., after transferring the spares), as previously introduced in \cite{han2025indirect}.

\subsubsection{Periodic Transition Matrix of Parking Orbits }
In this subsection, we derive the periodic transition structure of the parking orbit Markov chain, following the formulation for the constellation orbits in Sec.~\ref{sec:transition_const}. Specifically, we construct the transitions between the weighted distributions \(\tilde{\pi}^{\text{io}_\text{p}}_k\) and \(\tilde{\pi}^{\text{lt}_\text{p}}_{k,m}\) using the same time-resolved approach. In terms of notation, \((\cdot)_\text{d}\) is replaced with \((\cdot)_\text{i}\) (direct to indirect) and \((\cdot)_\text{c}\) with \((\cdot)_\text{p}\) (constellation to parking) for all related variables.

\paragraph{Transition into IO Phase} 
We assume that the parking spare inventory follows a single-source \((r,q,\tau)\) policy with no failures of parking spares. The update accounts for inflow from both the IO phase and the LT phase at \( m = m_\text{i} \), where indirect replenishment occurs with probability \( \rho_\text{i} \). The resulting IO-phase distribution at time \( k+1 \) is given by\footnote{To account for spare-satellite failures, multiply by the appropriate failure matrix at each time step, as in the constellation orbit analysis.}
\begin{equation} \label{eq:to_io_park}
\tilde{\pi}^{\text{io}_\text{p}}_{k+1} = 
\begin{cases}
\tilde{\pi}^{\text{io}_\text{p}}_k + \rho_\text{i} P_{q_\text{p}} \tilde{\pi}^{\text{lt}_\text{p}}_{k,m_\text{i}} ,
& k = 1,\dots,k_\text{p} - 1, \\
C^+_{r_\text{p}} \left[ 
P_{f_\text{p}} \left( \tilde{\pi}^{\text{io}_\text{p}}_{k_\text{p}} + \rho_\text{i} P_{q_\text{p}} \tilde{\pi}^{\text{lt}_\text{p}}_{k_\text{p},m_\text{i}} \right) 
\right],
& k = k_\text{p}.
\end{cases}
\end{equation}
Note that \( P_{f_\text{p}} \) is applied only at the alignment step, \( k = k_\text{p} \), when co-planar transfers to the constellation orbits are possible.

\paragraph{Transition into and within the LT Phase}
The transition into the LT phase with \( m = 0 \) at time step \( k+1 \) occurs only at \( k = k_\text{p} \) under the \((r,q,\tau)\) policy, i.e., orders are reviewed only at \( k = k_\text{p} \). A new order is initiated when the parking stock is below the reorder point, and the resulting distribution includes inflow from both the IO phase and the LT phase at \( m = m_\text{i} \). The update equation is
\begin{equation}
\tilde{\pi}^{\text{lt}_\text{p}}_{k+1,\,0}
=
\begin{cases}
  0,
    &  k = 1,\dots,k_\text{p} - 1, \\
  C^-_{r_\text{p}} \left[ P_{f_\text{p}} \left( \tilde{\pi}^{\text{io}_\text{p}}_k + \rho_\text{i} P_{q_\text{p}} \tilde{\pi}^{\text{lt}_\text{p}}_{k,m_\text{i}} \right) \right],
    &  k = k_\text{p}.
\end{cases}
\end{equation}

Next, the LT-to-LT transition for \( m < m_\text{i} \) simply advances the delay counter. At the RAAN alignment step, parking spares may be transferred to the constellation orbits, and the resulting depletion is applied via \( P_{f_\text{p}} \). The update equation is
\begin{equation}
\tilde{\pi}^{\text{lt}_\text{p}}_{k+1,\,m+1}
=
\begin{cases}
  \tilde{\pi}^{\text{lt}_\text{p}}_{k,m}, 
    &   k = 1,\dots,k_\text{p} - 1, \\
  P_{f_\text{p}}\,\tilde{\pi}^{\text{lt}_\text{p}}_{k,m}, 
    &   k = k_\text{p},
\end{cases}
\end{equation}
for \( m < m_\text{i} \).

As in the previous case, when \( m = m_\text{i} \), the update accounts for both advancement from the previous delay step and the probability of remaining at the current step without replenishment. The resulting update is given by
\begin{equation}
\tilde{\pi}^{\text{lt}_\text{p}}_{k+1,\,m_\text{i}}
=
\begin{cases}
  \tilde{\pi}^{\text{lt}_\text{p}}_{k,\,m_\text{i}-1}
  + (1-\rho_\text{i})\,\tilde{\pi}^{\text{lt}_\text{p}}_{k,\,m_\text{i}},
    & k = 1,\dots,k_\text{p}-1, \\
  P_{f_\text{p}} \left( \,\tilde{\pi}^{\text{lt}_\text{p}}_{k,\,m_\text{i}-1}
  + (1-\rho_\text{i})\,\,\tilde{\pi}^{\text{lt}_\text{p}}_{k,\,m_\text{i}} \right),
    & k = k_\text{p}.
\end{cases}
\end{equation}

\paragraph{Full Transition between IO and LT Phase}
Based on the update equations, we construct the full transition matrix by concatenating all sub-state distributions into a single vector:
\begin{equation}
    \boldsymbol{\pi}^\text{p}_k = \left(\tilde{\pi}^{\text{io}_\text{p}}_{k},\
    \tilde{\pi}^{\text{lt}_\text{p}}_{k,\,0},\,\dots,\,
    \tilde{\pi}^{\text{lt}_\text{p}}_{k,\,m_\text{i}}\right),
\end{equation}
where \( \boldsymbol{\mathbbm{1}} \cdot \boldsymbol{\pi}^\text{p}_k = 1 \).
The transition matrix for the case without RAAN alignment ($k \neq k_\text{p}$) is defined as:
\begin{equation} \label{eq:P_p_reg_mat}
    \mathbf{P}_\text{p,reg} = 
    \begin{bmatrix}
        I & 0  & \cdots & 0 & \rho_\text{i} P_{q_\text{p}} \\
        0 & 0  & \cdots  & 0 & 0 \\
        0 & I & \cdots & 0 & 0 \\
        \vdots & \vdots & \ddots  & \vdots & \vdots \\
        0 & 0 &  \cdots &  I & (1-\rho_\text{i})I
    \end{bmatrix}
\end{equation}
When RAAN alignment occurs (\(k = k_\text{p}\)), the transition matrix becomes:
\begin{equation} \label{eq:P_p_raan_mat}
    \mathbf{P}_{\text{p,raan}} = 
    \begin{bmatrix}
        C^+_{r_\text{p}} P_{f_\text{p}} & 0  & \cdots  & 0 & \rho_\text{i} C^+_{r_\text{p}} P_{f_\text{p}} P_{q_\text{p}} \\
        C^-_{r_\text{p}} P_{f_\text{p}} & 0  & \cdots  & 0 & \rho_\text{i} C^-_{r_\text{p}} P_{f_\text{p}} P_{q_\text{p}} \\
        0 & P_{f_\text{p}} & \cdots  & 0 & 0 \\
        \vdots & \vdots & \ddots  & \vdots & \vdots \\
        0 & 0 &  \cdots &  P_{f_\text{p}} & (1-\rho_\text{i})P_{f_\text{p}}
    \end{bmatrix}
\end{equation}
\textcolor{black}{Here, \(\mathbf{P}_\text{p,reg}, \mathbf{P}_{\text{p,raan}} \in \mathbb{R}^{(m_\text{i}+2)(N_{\text{sat}_\text{p}}+1) \times (m_\text{i}+2)(N_{\text{sat}_\text{p}}+1)}\).}

Using this formulation, the update rule for \(\boldsymbol{\pi}^\text{p}_k\) becomes:
\begin{equation} \label{eq:update_eq_parking}
    \boldsymbol{\pi}^\text{p}_{k+1} = 
    \begin{cases}
        \mathbf{P}_\text{p,reg}\,\boldsymbol{\pi}^{\text{p}}_{k}, 
        & \text{if } k = 1,\dots,k_\text{p}-1, \\
        \mathbf{P}_\text{p,raan}\,\boldsymbol{\pi}^{\text{p}}_{k}, 
        & \text{if } k = k_\text{p}.
    \end{cases}
\end{equation}
Finally, the periodic steady-state condition is expressed as:
\begin{equation} \label{eq:stationary_parking}
    \boldsymbol{\pi}^\text{p}_{1} = \mathbf{P}_\text{p,raan}\left( \mathbf{P}_\text{p,reg}  \right)^{k_\text{p}-1} \boldsymbol{\pi}^\text{p}_{1},
\end{equation}
\textcolor{black}{which serves as a fixed-point equation to determine the long-run state distribution of the parking-orbit inventory over a full RAAN-alignment cycle.}

\begin{remark}
As explained earlier, the time-resolved method is systematic and intuitive but results in a state-space dimension that scales with \( m_\text{i} \), which can become large when the time step \( \tau_\text{mc} \) is small or the processing time \( \tau_\text{lv,i} \) is long. In Sec.~\ref{sec:approx_method}, we introduce a simple approximation method to reduce this dimensionality.
\end{remark}

%
\subsubsection{Average Distribution over the Replenishment Cycle}
\textcolor{black}{The periodic stationary distribution \( \boldsymbol{\pi}^\text{p}_{k} \) at each step \( k \) of the RAAN-alignment cycle is obtained by solving Eq.~\eqref{eq:stationary_parking}.} Let \( \mathbf{I}_\text{p} \) denote the row-block matrix
\begin{equation}
    \mathbf{I}_\text{p} = \begin{bmatrix} I_{N_{\text{sat}_\text{p}} + 1} & \cdots & I_{N_{\text{sat}_\text{p}} + 1} \end{bmatrix},
\end{equation}
with \( m_\text{i} + 2 \) identity blocks corresponding to the IO phase and all \( m_\text{i} + 1 \) LT sub-phases. Then, the average distribution at step \( k \) is
\begin{equation}
    \pi^{\text{rc}_\text{p}}_k = \mathbf{I}_\text{p} \boldsymbol{\pi}^{\text{p}}_k = 
    \tilde{\pi}^{\text{io}_\text{p}}_{k} + \sum_{m=0}^{m_\text{i}} \tilde{\pi}^{\text{lt}_\text{p}}_{k,m},
\end{equation}
and the average distribution over the replenishment cycle is
\begin{equation}
    \pi^{\text{rc}_\text{p}} = \frac{1}{k_\text{p}} \sum_{k=1}^{k_\text{p}} \pi^{\text{rc}_\text{p}}_k.
\end{equation}

%
\subsubsection{Computation of Average Replenishment Period}

This approach follows the same structure as the computation for constellation orbits. Define
\begin{equation}
    c^{\text{io}_\text{p}} = \sum_{k=1}^{k_\text{p}} \mathbbm{1} \cdot \tilde{\pi}^{\text{io}_\text{p}}_{k}, \qquad
    c^{\text{lt}_\text{p}} = \sum_{k=1}^{k_\text{p}} \sum_{m=0}^{m_\text{i}} \mathbbm{1} \cdot \tilde{\pi}^{\text{lt}_\text{p}}_{k,m}.
\end{equation}
From the lead-time model, the average LT duration is
\begin{equation}
    \tau_{\text{lt}_\text{p}} = m_\text{i} \tau_\text{mc} + \mu_\text{lv,i},
\end{equation}
and the corresponding IO duration is
\begin{equation}
    \tau_{\text{io}_\text{p}} = \frac{c^{\text{io}_\text{p}}}{c^{\text{lt}_\text{p}}} \tau_{\text{lt}_\text{p}}.
\end{equation}
Thus, the average replenishment cycle time is
\begin{equation}
    \tau_{\text{rc}_\text{p}} = \tau_{\text{lt}_\text{p}} + \tau_{\text{io}_\text{p}}.
\end{equation}

%

\subsubsection{Parking Spares Availability Distribution} \label{sec:park_avail}

To compute the parking-spares availability at the RAAN-contact event \( E \), we require the inventory distribution immediately after indirect replenishment (\( P_{q_\text{p}} \)) but before any transfer to the constellation orbits (\( P_{f_\text{p}} \)). As in the demand distribution analysis at event \( E \) in Sec.~\ref{sec:dmd_dist}, we apply the regular transition \( \mathbf{P}_\text{p,reg} \) to obtain this pre-transfer state. Accordingly, the average distribution conditioned on event \( E \) is
\begin{equation} \label{eq:pi_rc_E}
    \pi^{\text{rc}_\text{p}|E} = \mathbf{I}_\text{p}\,
        \mathbf{P}_\text{p,reg}\,\boldsymbol{\pi}^{\text{p}}_{k_\text{p}},
\end{equation}
and the availability probability defined in Eq.~\eqref{eq:parking_avail} becomes
\begin{equation} \label{eq:kappa_eq}
    \kappa_j = \sum_{k=j}^{N_{\text{sat}_\text{p}}}\pi^{\text{rc}_\text{p}|E} (X_\text{p} = k),
\end{equation}
where \( \kappa_j \) is the probability that at least \( j \) parking spares are available at the RAAN-alignment moment.

\subsection{Flow of Hybrid Strategy Analysis}
This paper adopts the same fixed-point solution approach introduced in \cite{han2025indirect}. The coupled analysis of in-plane and parking orbits is separated into two parts, as described in the preceding subsections. These components are solved iteratively, passing updated values of \( \chi \) and \( \kappa \) between them, starting from full parking availability (\( \kappa = \mathbbm{1} \)), until a consistent solution is reached.

The most computationally intensive steps involve solving the linear systems in Eq.~\eqref{eq:stationary_constel} and Eq.~\eqref{eq:stationary_parking}. However, for practical parameter settings (summarized in Sec.~\ref{sec5}), the solution can be obtained in under a second without approximation, and in just a few milliseconds when the approximation method is used, even on a standard desktop computer.

The complete procedure for the hybrid resupply strategy is summarized in Table~\ref{alg:cap}. Note that the parking orbit analysis step (from $P_{f_\text{p}}$ to $\kappa$) can be replaced by the event-based method proposed in \cite{han2025indirect}.
\begin{algorithm}[H]
\caption{Fixed Point Iteration for the Hybrid Strategy}\label{alg:cap}
\begin{algorithmic}
\Require Constellation Configuration, Probability Model
\State Initialize $\tau_\text{c}$, $\tau_\text{p}$, $P_{f}$, $\rho_\text{d}$, $\rho_\text{i}$, $C_{r_\text{c,d}}^\pm$, $P_{q_\text{c,d}}$, and $P_{q_\text{p}}$ based on Sec.~\ref{sec2}
\State $\kappa^{(0)} \gets \mathbbm{1}$ and $n \gets 0$
\Repeat
    \State $P_{q_\text{c,i}}\gets$ Eq.~\eqref{eq:park_avail_mat} using $\kappa^{(n)}$
    \State $\mathbf{P}_{\text{c},\text{reg}} \gets$ Eq.~\eqref{eq:P_c_reg_mat} and $\mathbf{P}_{\text{c},\text{raan}}\gets$ Eq.~\eqref{eq:P_c_raan_mat}
    \State $\boldsymbol{\pi}^\text{c}_{1} \gets$ Eq.~\eqref{eq:stationary_constel} and $\boldsymbol{\pi}^\text{c}_{k} \gets$ Eq.~\eqref{eq:update_eq_constel} for $k > 1$
    \State $\pi^{\text{rc}_\text{c}|E} \gets$ Eq.~\eqref{eq:pi_r_c_i} and  $\chi \gets$ Eq.~\eqref{eq:inplane_demand_prob}
    \State $P_{f_\text{p}} \gets $ Eq.~\eqref{eq:fail_matrix_park}
    \State $\mathbf{P}_{\text{p},\text{reg}} \gets$ Eq.~\eqref{eq:P_p_reg_mat} and $\mathbf{P}_{\text{p},\text{raan}}\gets$ Eq.~\eqref{eq:P_p_raan_mat}
    \State $\boldsymbol{\pi}^\text{p}_{1} \gets$ Eq.~\eqref{eq:stationary_parking} and $\boldsymbol{\pi}^\text{p}_{k} \gets$ Eq.~\eqref{eq:update_eq_parking} for $k > 1$
    \State $\pi^{\text{rc}_\text{p}|E} \gets$ Eq.~\eqref{eq:pi_rc_E} and $\kappa^{(n+1)} \gets$ Eq.~\eqref{eq:kappa_eq}
    \State $n \gets n + 1$
\Until{$\|\kappa^{(n)}-\kappa^{(n-1)}\| \le \varepsilon$ or $n \ge n^{\max}$}
\end{algorithmic}
\end{algorithm}
\begin{remark}
    The optimal numerical method for solving the stationary distribution depends on the problem scale. For the large, sparse matrices of the full-scale analysis, an iterative eigensolver (e.g., a Krylov subspace method \cite{greenbaum1997iterative}) with matrix-vector products is more efficient than explicitly computing the matrix power as in Eq.~\eqref{eq:stationary_constel}. For the smaller system resulting from the approximation in Sec.~\ref{sec:approx_method}, it is faster to form the full transition matrix and solve $(\mathbf{I} - \mathbf{P}_{{(\cdot)},\text{trn}})\boldsymbol{\pi}_1^{(\cdot)} = \mathbf{0}$ together with the normalization condition $\mathbbm{1}\cdot\boldsymbol{\pi}_1^{(\cdot)}=1$.
\end{remark}

\subsection{Approximate Modeling for Efficient Computation} \label{sec:approx_method}
To reduce the dimensionality of the linear systems in Eqs.~\eqref{eq:stationary_parking} and~\eqref{eq:stationary_constel}, we introduce a simple approximation technique. While more sophisticated methods based on exact distributions are possible, they increase the overall complexity. In contrast, the proposed method is easy to apply and sufficiently accurate across various numerical tests.

\subsubsection{Reduction of the dimension of in-plane stock}
Because dual-channel replenishment makes deep stock outages extremely unlikely, we truncate the state space below a conservative lower bound.  Let
\[
\tau_\text{slow} = \max\{\tau_\text{c},\;\mu_{\text{lv,d}} + \tau_{\text{lv,d}}\}
\quad\text{and}\quad
M_0 = \min\{N_{\text{sat}_\text{c}},\;r_\text{c,i},\;r_\text{c,d}\}.
\]
Over a horizon of
\(\;n=\tau_\text{slow}/\tau_\text{mc}\), the expected failures are
\[
\mathbb{E}[F] = n\,M_0\,\lambda_\text{sat},\quad {\sigma}[F] = \sqrt{n\,M_0\,\lambda_\text{sat}}.
\]
We then set the effective lower bound
\[
M_{\text{sat}_\text{c}} = \max\bigl[ 
\lfloor M_0 - \mathbb{E}[F] \;-\; k_\sigma{\sigma}[F]\rfloor,\, 0\bigr],
\]
which accounts for both the expected failures and $k_\sigma$ standard deviations as a safety margin.

We then restrict the in-plane stock state to  
\(\;X_\text{c}\in \{ M_{\text{sat}_\text{c}},\dots, N_{\text{sat}_\text{c}} \}\),  
and renormalize every transition matrix (e.g.\ \(P_f\), \(P_{q_\text{c,i}}\), …) so each column sums to one.

\subsubsection{Approximation of the Launch Delay Model}
The dimension of the linear system is proportional to the number of delay stages \( m_{(\cdot)} \), which corresponds to the discretized launch delay \( \tau_{\text{lv},(\cdot)} \) with step size \( \tau_\text{mc} \). For simplicity, we explain the in-plane case; the same reduction applies to the parking-orbit model. To reduce the system size, we approximate the original \( m \)-step delay chain using a smaller number of equally spaced stages \( s \), each modeled as a first-order lag process.

Let \( m_{\ell} \) denote the number of discrete time steps in the \( \ell \)-th reduced stage, representing the delay between \( \pi^{\text{lt}_\text{c}}_{k,\,\ell-1} \) and \( \pi^{\text{lt}_\text{c}}_{k,\,\ell} \), such that $\sum_{\ell=1}^{s} m_\ell = m_\text{d}$.
For a user-defined number of reduced stages $ s \leq m_{\text{d}} $, the stage durations are chosen to distribute the total delay as evenly as possible: the first $ r = m_{\text{d}} - s \cdot \lfloor m_{\text{d}}/s \rfloor $ stages are assigned $ \lfloor m_{\text{d}}/s \rfloor + 1 $ steps each, and each of the remaining $ s - r $ stages receive $ \lfloor m_{\text{d}}/s \rfloor $ steps.

With this new definition, the reduced state vector becomes
\[
\boldsymbol{\pi}^\text{c}_k = \left(\tilde{\pi}^{\text{io}_\text{c}}_{k},\
\tilde{\pi}^{\text{lt}_\text{c}}_{k,\,0},\,\dots,\,
\tilde{\pi}^{\text{lt}_\text{c}}_{k,\,s}\right),
\]
reducing the state size from \( (m_\text{d}+2)(N_{\text{sat}_\text{c}}+1) \) to \( (s+2)(N_{\text{sat}_\text{c}}+1) \).

Applying the delay-stage reduction to Eq.~\eqref{eq:P_c_reg_mat} gives the approximate regular-step update. In the full model, the LT phase explicitly tracks every one-step delay state. In the approximate model, \(m_\ell\) consecutive full-model LT sub-phases are grouped into one reduced stage: a fraction \(1/m_\ell\) advances to the next reduced stage after the accumulated stock transition over those \(m_\ell\) steps, while the remaining fraction \(1-1/m_\ell\) stays in the current reduced stage. Therefore, for \(k\neq k_\text{c}\),
\begin{equation}
\tilde{\pi}^{\text{lt}_\text{c}}_{k+1,\,\ell} = 
\begin{cases}
    C_{r_\text{c,d}}^- P_f \left( \tilde{\pi}^{\text{io}_\text{c}}_k + \rho_\text{d} P_{q_\text{c,d}} \tilde{\pi}^{\text{lt}_\text{c}}_{k,s}\right) + \left( 1 - \dfrac{1}{m_1}\right) \tilde{\pi}^{\text{lt}_\text{c}}_{k,0}, & \ell = 0, \\
    \dfrac{1}{m_{\ell}} \left(P_f\right)^{m_\ell} \tilde{\pi}^{\text{lt}_\text{c}}_{k,\,\ell -1}+\left( 1 - \dfrac{1}{m_{\ell+1}}\right)\tilde{\pi}^{\text{lt}_\text{c}}_{k,\,\ell},
    & \ell = 1,\dots,s-1, \\
    \dfrac{1}{m_s} \left(P_f\right)^{m_s}\tilde{\pi}^{\text{lt}_\text{c}}_{k,\,s-1}+\left( 1 - \rho_\text{d} \right) P_f \tilde{\pi}^{\text{lt}_\text{c}}_{k,\,s},
    & \ell = s.
\end{cases}
\end{equation}
The same delay-stage reduction is applied to Eq.~\eqref{eq:P_c_raan_mat} to obtain the approximate RAAN-step matrix. Thus, the approximation preserves the regular/RAAN event sequence defined in the full matrices and only changes how the LT-delay chain is represented. The parking-orbit version is obtained by applying the same reduction to the corresponding parking-orbit transition matrices.

As a concrete example, let \( \tau_\text{mc} = 1 \), \( \tau_\text{lv,d} = 8 \), and \( s = 3 \). Then, the reduced stage durations become \( (m_1, m_2, m_3) = (3, 3, 2) \), and the approximated transition matrix is:
\[
\mathbf{P}^{\mathrm{approx}}_{\text{c},\text{reg}} = \begin{bmatrix}
    C^+_{r_\text{c,d}}P_{f} & 0   & 0 & 0 & \rho_\text{d} C^+_{r_\text{c,d}} P_{f} P_{q_\text{c,d}}\\
    C^-_{r_\text{c,d}}P_{f} & \frac{2}{3}I   & 0& 0 & \rho_\text{d} C^-_{r_\text{c,d}} P_{f} P_{q_\text{c,d}} \\
    0 & \frac{1}{3}P_{f}^3 & \frac{2}{3}I  & 0 & 0 \\
    0 & 0 &  \frac{1}{3}P_{f}^3 & \frac{1}{2}I   & 0 \\
    0 & 0 &  0  & \frac{1}{2}P_{f}^2 & (1-\rho_\text{d})P_{f} \\
\end{bmatrix}
\]
Each column of this matrix sums to one, preserving the Markov property under the reduced-stage approximation.

\section{Performance Evaluation of Spare Management Policy} \label{sec4}

Once the stationary distributions $\pi^{\text{rc}_{(\cdot)}}$ are obtained, we can assess key performance indicators of the replenishment strategy. Two principal metrics are considered: operational cost and resilience, which typically exhibit a trade-off relationship. For consistency and comparability, we extend the cost and constraint modeling framework previously introduced in \cite{han2025direct, han2025indirect}.

\subsection{Cost Model for Hybrid Strategy}

The overall expected cost per unit time is composed of four terms:
\begin{equation}
    C_\text{total} = C_\text{build} + C_\text{hold} + C_\text{launch} + C_\text{trans},
\end{equation}
where each component is defined as follows:

\paragraph{Manufacturing cost}
\begin{equation}
    C_\text{build} = c_\text{build} \left(
        \frac{N_{\text{orbit}_\text{c}} \, q_\text{c,d}}{\tau_{\text{rc}_\text{c}}}
        + \frac{N_{\text{orbit}_\text{p}} \, q_\text{c,i} \, q_\text{p}}{\tau_{\text{rc}_\text{p}}}
    \right)
\end{equation}
This term represents the expected manufacturing cost per unit time, with production rates determined by the average replenishment frequency for each channel. Here, \( c_\text{build} \) is the unit cost to manufacture a spare satellite.

\paragraph{Holding cost}
\begin{equation}
\begin{aligned}
    C_\text{hold} &= 
    c_{\text{hold}_\text{c}} N_{\text{orbit}_\text{c}} 
    \sum_{x=\bar{N}_\text{sat}+1}^{N_{\text{sat}_\text{c}}} (x - \bar{N}_\text{sat}) \ \pi^{\text{rc}_\text{c}}(X_\text{c}=x) \\
    & + c_{\text{hold}_\text{p}} N_{\text{orbit}_\text{p}} q_\text{c,i}  \sum_{x=0}^{N_{\text{sat}_\text{p}}} 
    x \ \pi^{\text{rc}_\text{p}}(X_\text{p}=x).
\end{aligned}
\end{equation}
This penalizes excess inventory in orbit. For constellation orbits, only the surplus above $\bar{N}_\text{sat}$ is treated as spare, while all parking stock is considered spare by design. Here, \( c_{\text{hold}_\text{c}} \) and \( c_{\text{hold}_\text{p}} \) are per-satellite holding costs for constellation and parking orbits.

\paragraph{Launch cost}
\begin{equation}
\begin{aligned}
    C_\text{launch} =
    \frac{N_{\text{orbit}_\text{c}}}{\tau_{\text{rc}_\text{c}}} c_{\text{full,d}}
    + \frac{N_{\text{orbit}_\text{p}}}{\tau_{\text{rc}_\text{p}}}c_{\text{full,i}}
\end{aligned}
\end{equation}
The total launch cost is the sum of the expected costs for dedicated launches in the direct and indirect channels, where \( c_{\text{full,d}} \) and \( c_{\text{full,i}} \) are the respective costs per launch vehicle. This model considers only dedicated launches, as rideshare is often impractical for on-demand replenishment of a specific constellation plane; therefore, the rideshare cost modeling from our previous work~\cite{han2025indirect} is not included.

\paragraph{Orbital transfer cost}
\begin{equation}
\begin{aligned}
    C_\text{trans} = 
    \frac{N_{\text{orbit}_\text{c}}}{\tau_\text{c} q_\text{c,i}} 
    (c_\text{fuel} \, m_\text{fuel,i} + c_\text{trans}) 
    \sum_{x=0}^{N_{\text{sat}_\text{c}}} x
    \left[
    \pi^{q_\text{c,i}}(X_\text{c}=x)
    -
    \pi^{r_\text{c,i}}(X_\text{c}=x)
    \right].
\end{aligned}
\end{equation}
Here, $c_\text{fuel}$ is the unit fuel cost, $c_\text{trans}$ covers non-fuel transfer expenses, and $m_\text{fuel,i}$ is computed from $\Delta v$ and dry mass, \(m_{\rm dry,i}=q_\text{c,i}\,m_{\rm sat}+m_{\rm bus}\),  using Eqs.~\eqref{eq:delta_v} and \eqref{eq:fuel_mass}. 

\subsection{Resilience Metric}
A proper resilience metric should capture both {agility} (how quickly the system recovers to nominal capacity) and {robustness} (the depth of degradation when operating below nominal) \cite{dod_resilience_2011, Ron_resilience}. In practice, this is often quantified by the area under the nominal capacity line lost over time, which corresponds to the time-integrated shortage.

In our discrete-time Markov model, this becomes the expected shortage in constellation orbits:
\begin{equation} \label{eq:shortage_plane}
    S_\text{c} = \sum_{i=0}^{\bar{N}_\text{sat}} (\bar{N}_\text{sat} - i) \cdot \pi^{\text{rc}_\text{c}}(i),
\end{equation}
which weighs the deficit $(\bar{N}_\text{sat} - i)$ by the steady state probability of being in state $i$.

In parking orbits, where no nominal threshold is defined, we instead use the probability of being out-of-stock:
\begin{equation} \label{eq:res_park}
    \mathbb{P}(X_\text{p} = 0) = \pi^{\text{rc}_\text{p}}(0),
\end{equation}
which represents the long-run fraction of time the parking inventory is empty.

\subsection{Optimization Problem}
We seek to minimize the total expected cost while satisfying resilience and launch capacity constraints:
\begin{equation} \label{eq:opt_formulation}
\begin{aligned}
    \min_{x} \quad & C_\text{total} \\
    \text{s.t.} \quad 
        & g_1 = S_\text{c} - \varepsilon_1 \leq 0 \quad \text{(constellation resilience)} \\
        & g_2 = \mathbb{P}(X_\text{p} = 0) - \varepsilon_2 \leq 0 \quad \text{(parking availability)} \\
        & g_3 = m_\text{total,d} - m_\text{payload,d} \leq 0 \quad \text{(direct LV capacity)} \\
        & g_4 = m_\text{total,i} - m_\text{payload,i} \leq 0 \quad \text{(indirect LV capacity)} \\
        & q_\text{c,d}, q_\text{c,i}, r_\text{c,d}, r_\text{c,i}, q_\text{p}, r_\text{p}, N_{\text{orbit}_\text{p}} \in \mathbb{Z}^+
\end{aligned}
\end{equation}
Here, the decision vector $x$ includes batch sizes $(q_\text{c,d}, q_\text{c,i}, q_\text{p})$, reorder points $(r_\text{c,d}, r_\text{c,i}, r_\text{p})$, the number of parking orbits $(N_{\text{orbit}_\text{p}})$, and the altitude of the parking orbits ($h_\text{p}$). The thresholds \( \varepsilon_1 \) and \( \varepsilon_2 \) define the maximum allowable shortage and out-of-stock probability, respectively. Constraints \( g_3 \) and \( g_4 \) ensure that the total mass per launch does not exceed the maximum payload capacity for the direct ($m_\text{payload,d}$) and indirect ($m_\text{payload,i}$) launch vehicles. The total mass for each channel is defined as $ m_\text{total,d} = m_\text{sat} q_\text{c,d}$ and $ m_\text{total,i} = (m_\text{fuel,i} + m_\text{dry,i})q_\text{p}$. For simplicity, a constant maximum payload capacity is assumed for each vehicle, regardless of its destination orbit.

\section{Numerical Validation of the Analysis Method} \label{sec5}
\subsection{Numerical Validation Setup}
When formulating the model, we assumed that stock levels in both constellation and parking orbits are i.i.d. Since directly comparing the full distribution $\pi^{\text{rc}_{(\cdot)}}$ between analysis and simulation is impractical, we instead use representative performance metrics. Inspired by the validation approach in \cite{jakob2019optimal}, we compare the mean stock levels of in-plane and parking orbits, the expected shortage of the constellation orbit (Eq.~\eqref{eq:shortage_plane}), and the out-of-stock probability of the parking orbit (Eq.~\eqref{eq:res_park}). The mean stock level is computed as  
\begin{equation}
M_{(\cdot)} = \sum_{i=0}^{N_{\text{sat}_{(\cdot)}}} i \cdot \pi^{\text{rc}_{(\cdot)}}(X_{(\cdot)}=i).
\end{equation}

We generate 100 test cases using Latin hypercube sampling over the parameter ranges in Table~\ref{tab:trade_space_lhs}, with fixed parameters listed in Table~\ref{tab:fixed_sim_para}. \textcolor{black}{The fixed parameters represent a hypothetical mid-inclination LEO broadband mega-constellation at a tractable but representative scale for validation.} Each case is simulated for 20 years, repeated 100 times, and averaged to obtain the statistics. For $M_\text{p}$ and $\mathbb{P}(X_\text{p}=0)$, 16 cases with $r_\text{c,d} > r_\text{c,i}+4$ were excluded, since indirect usage was too small for a 20-year simulation to capture the very low-probability transitions reliably\textcolor{black}{, reflecting a finite-horizon Monte Carlo limitation for rare events}.

\begin{table}[hbt!]
\centering
\footnotesize
\setlength{\tabcolsep}{5pt}
\renewcommand{\arraystretch}{1.1}
\caption{Fixed simulation parameters}
\label{tab:fixed_sim_para}
\begin{tabularx}{\linewidth}{@{}>{\raggedright\arraybackslash}X c c c@{}}
\toprule
Parameter & Notation & Value & Unit \\

\midrule
Markov time step
  & $\tau_\text{mc}$ 
  & $0.5$ 
  & days \\

Constellation orbit altitude ($a_\text{c} - R_\oplus$)
  & $h_\text{c}$ 
  & $1200$ 
  & km \\
  
Inclination of orbit planes
  & $i$ 
  & $50$ 
  & deg \\

Number of constellation orbits
  & $N_{\text{orbit}_\text{c}}$ 
  & $40$ 
  & orbits \\

Nominal satellites per plane
  & $\bar{N}_\text{sat}$ 
  & $40$ 
  & satellites \\

\bottomrule
\end{tabularx}
\end{table}

\begin{table}[hbt!]
\centering
\footnotesize
\setlength{\tabcolsep}{5pt}
\renewcommand{\arraystretch}{1.1}
\caption{Bounds of sampled simulation parameters}
\label{tab:trade_space_lhs}
\begin{tabularx}{\linewidth}{@{}>{\raggedright\arraybackslash}X c c c@{}}
\toprule
Parameter & Notation & Bounds & Unit \\
\midrule
Annual satellite failure rate & $\lambda_{\text{sat,yr}}$ & $[0.001,\,0.5]$ & fail/sat/yr \\
Direct LV order processing time & $\tau_{\text{lv,d}}$ & $[0,\,60]$ & days \\
Indirect LV order processing time & $\tau_{\text{lv,i}}$ & $[0,\,60]$ & days \\
Mean direct LV exponential lead time & $\mu_{\text{lv,d}}$ & $[5,\,60]$ & days \\
Mean indirect LV exponential lead time & $\mu_{\text{lv,i}}$ & $[5,\,60]$ & days \\
Order size, in-plane (indirect channel) & $q_{\text{c,i}}$ & $[1,\,10]$ & satellites \\
Order size, in-plane (direct channel) & $q_{\text{c,d}}$ & $[1,\,10]$ & satellites \\
Order size, parking-orbit batch & $q_{\text{p}}$ & $[1,\,20]$ & batches \\
Reorder point, in-plane (indirect channel) & $r_{\text{c,i}}$ & $[30,\,45]$ & satellites \\
Reorder point, in-plane (direct channel) & $r_{\text{c,d}}$ & $[30,\,45]$ & satellites \\
Reorder point, parking-orbit batch & $r_{\text{p}}$ & $[0,\,10]$ & batches \\
Parking-orbit altitude ($a_\text{p}-R_\oplus$) & $h_\text{p}$ & $[500,\,1100]$ & km \\
Number of parking orbits & $N_{\text{orbit}_\text{p}}$ & $[1,\,10]$ & orbits \\
\bottomrule
\end{tabularx}
\end{table}

\subsection{Numerical Validation Results} \label{sec:num_val}
As noted in our earlier work \cite{han2025indirect}, the proposed analysis can lose accuracy when the i.i.d. assumption for parking orbits breaks down, which occurs if parking spares remain out-of-stock for an extended period. This typically happens when the satellite failure rate exceeds the replenishment capacity of the launch vehicle, disrupting the conventional saw-tooth (i.e., even-distribution) stock profile. The heuristic threshold $\mathbb{P}(X_\text{p} = 0) < 1/(N_{\text{sat}_\text{p}}+1)$ is derived from this even-distribution assumption and is used here to select cases for comparison. 

The results in Table~\ref{tab:error_list} confirm the high fidelity of the analytical model, with the 95th percentile error below 1\% for all metrics. Discrepancies arise from the i.i.d. assumption in the model and statistical noise in finite-duration Monte Carlo runs. The larger relative errors for $M_\text{p}$ and $S_\text{c}$ reflect their small nominal values, which is a known limitation of relative error metrics. 

\begin{figure}[!h]
    \centering
    \includegraphics[width=.45\textwidth]{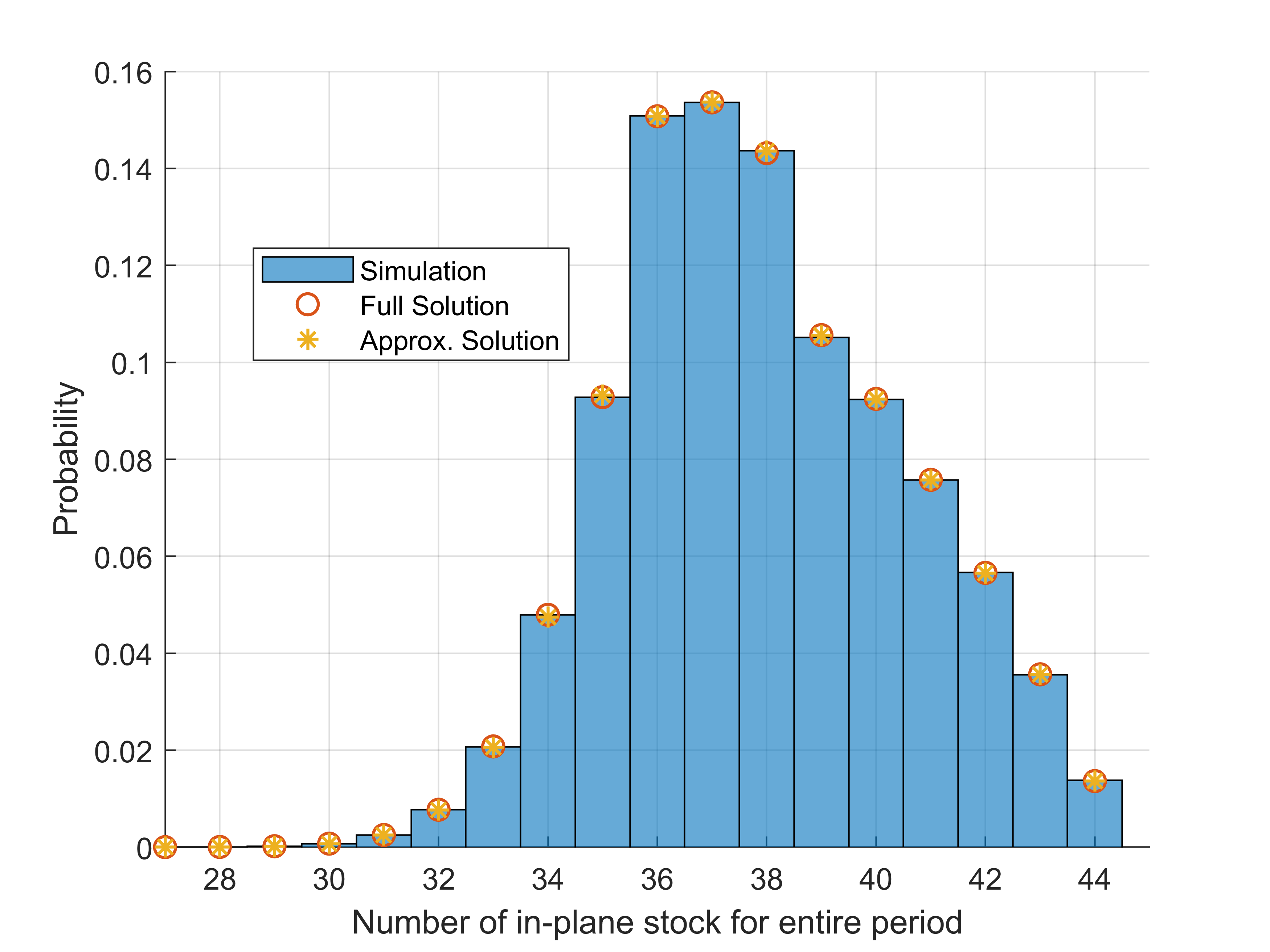}
    \includegraphics[width=.45\textwidth]{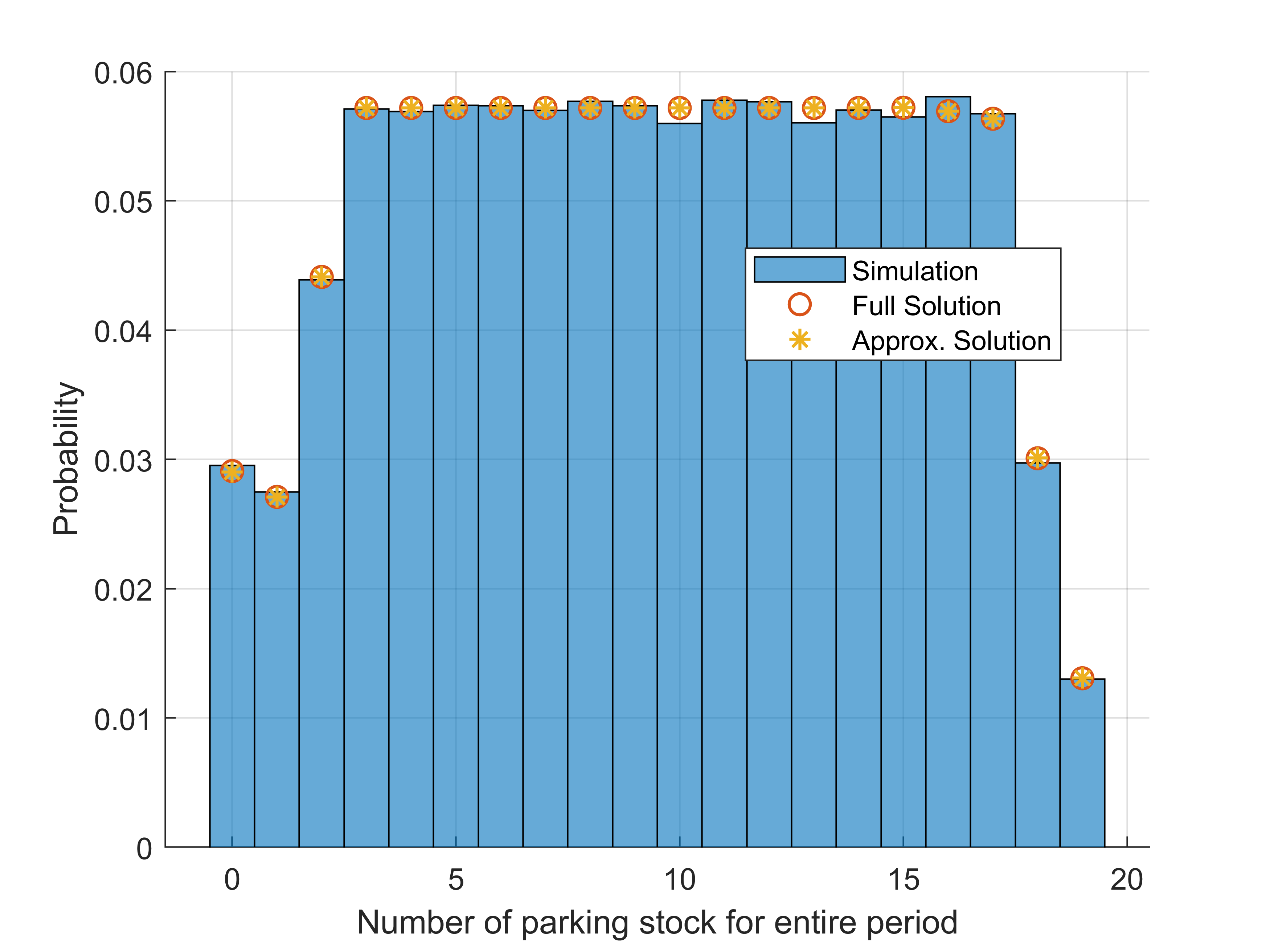}
    \caption{Comparison results for the representative case of large error: (a) $\pi^{\text{rc}_\text{c}}$; (b) $\pi^{\text{rc}_\text{p}}$}
    \label{fig:worst_case_result}
\end{figure}

A representative case near the maximum error is shown in Fig.~\ref{fig:worst_case_result}. The in-plane distribution matches almost perfectly, while the parking orbit distribution exhibits minor oscillations around the analytical points, a pattern consistent with Monte Carlo sampling noise; extending the simulation horizon reduces the discrepancy, confirming this interpretation.

In terms of computational cost, the full Monte Carlo simulation required a few hours per test case, whereas the analytical framework produced a solution in under a second with appropriate numerical techniques. The method is therefore both accurate and highly efficient, making it suitable for use as the inner loop in an optimization process.

\begin{table}[hbt!]
\centering
\caption{Error between proposed method and simulation results}
\label{tab:error_list}
\begin{tabular}{lc c c }
\toprule
Parameter & Mean & P95\\
\midrule
Relative error of $M_\text{c}$ 
  & $0.009 \ \%$ 
  & $0.023 \ \%$ \\

Relative error of $M_\text{p}$ 
  & $0.309 \ \%$ 
  & $0.776 \ \%$ \\
  
Relative error of $S_\text{c}$ 
  & $0.151 \ \%$ 
  & $0.545 \ \%$ \\

Absolute error of $\mathbb{P}(X_\text{p} = 0)$ 
  & $0.023 \text{ p.p}$ 
  & $0.092 \text{ p.p}$ \\

\bottomrule
\end{tabular}
\end{table}

\subsection{Effectiveness of the Approximate Solution} \label{sec:approx_analysis}
The approximation introduced in Sec.~\ref{sec:approx_method} is designed to reduce computational cost while maintaining accuracy, enabling its efficient use within the optimization process of Sec.~\ref{sec6}. The approximation parameters were determined using the 100 test cases from the previous validation. It was found that a safety margin based on a standard deviation of $k_\sigma = 3$ was sufficient to account for nearly all observed stock-level drops. Likewise, a sensitivity analysis of the stage-reduction parameter $s$ indicated that $s = \lfloor m_\text{lv,d}/3 \rceil$ offers a good balance between accuracy and computation time.

With these parameters, the effectiveness of the proposed technique was confirmed. The difference between the approximated and full solutions was less than $0.05\%$ for all performance metrics, demonstrating that the approximation accurately captures the primary system dynamics. Furthermore, the approximation offers a substantial computational gain: the average solution time per case decreases from one second for the full model to a few tens of milliseconds, yielding a 20-fold speed-up\footnote{The speed comparison was performed using the fastest available numerical solver for each respective model: an iterative eigensolver for the full model and a direct linear solve for the approximation.}.

\subsection{Comparison with Existing Method}
\subsubsection{Limiting Assumptions in Prior Hybrid Strategy Analysis}
The studies in \cite{jakob2019optimal,Kim2025_dual} model demand of in-plane orbit as a Poisson process and assume a uniform lead-time distribution from parking to constellation orbits. As discussed in our earlier work on the indirect strategy \cite{han2025indirect}, these assumptions can break down and fail to capture the multi-echelon nature of the problem in general.

\begin{figure}[!h]
    \centering
    \includegraphics[width=.6\textwidth]{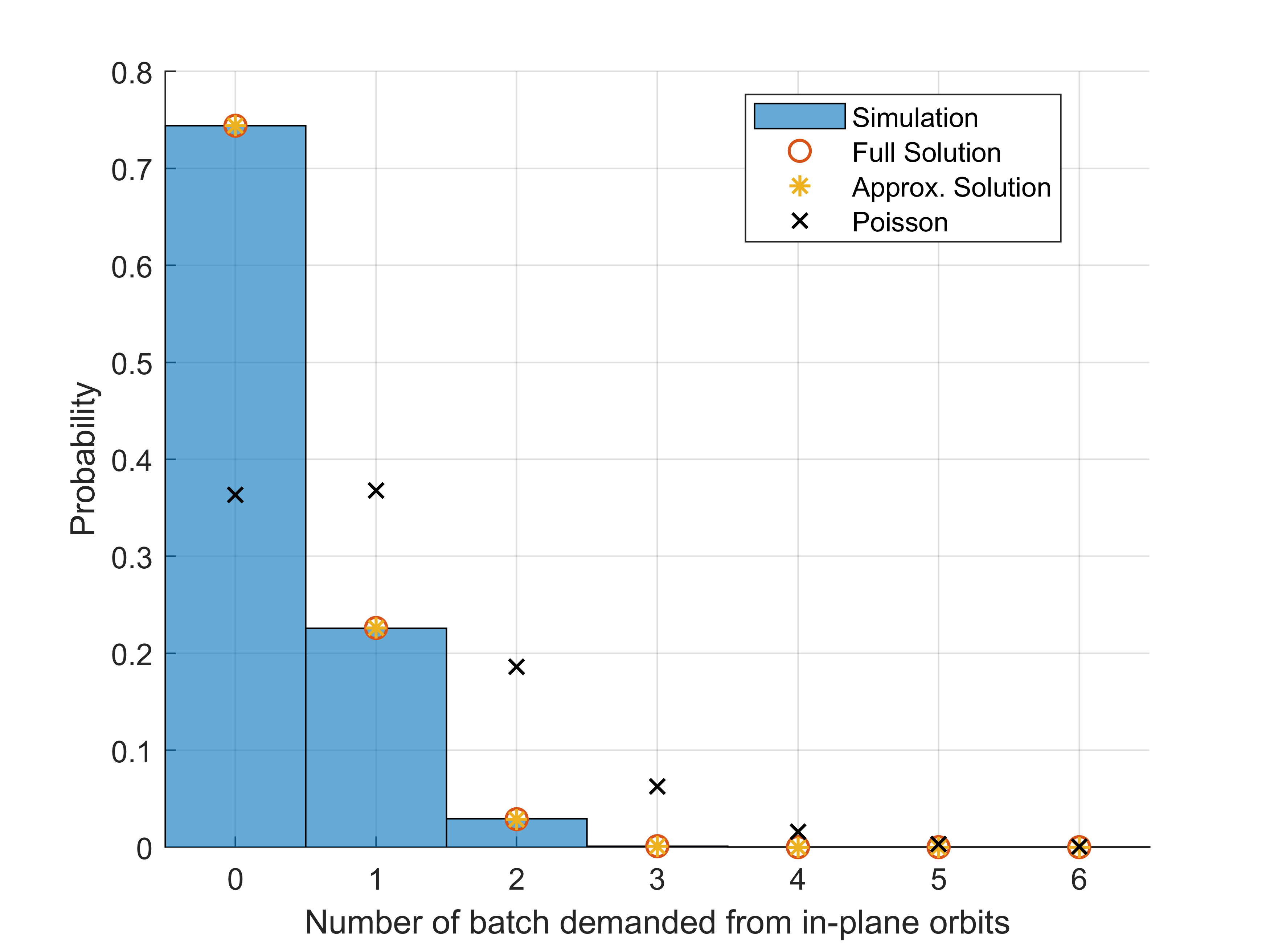}
    \caption{Parking-spare demand distribution from in-plane orbits for the representative case}
    \label{fig:demand_test}
\end{figure}

Figure~\ref{fig:demand_test} shows the demand distribution for parking spares from in-plane orbits, $\chi$. The black crosses indicate the corresponding Poisson model \cite{jakob2019optimal,Kim2025_dual}
\begin{equation*}
    \lambda_\text{demand} = \frac{N_{\text{orbit}_\text{c}} \bar{N}_{\text{sat}}}{q_\text{c,i} N_{\text{orbit}_\text{p}}} \lambda_\text{sat} \tau_\text{c}.
\end{equation*}
The actual distribution deviates significantly because (i) the direct channel influences parking demand, and (ii) demand depends on parking stock availability. In contrast, our proposed analysis (red circles and orange stars) matches the simulation closely.

\begin{figure}[!h]
    \centering
    \includegraphics[width=.6\textwidth]{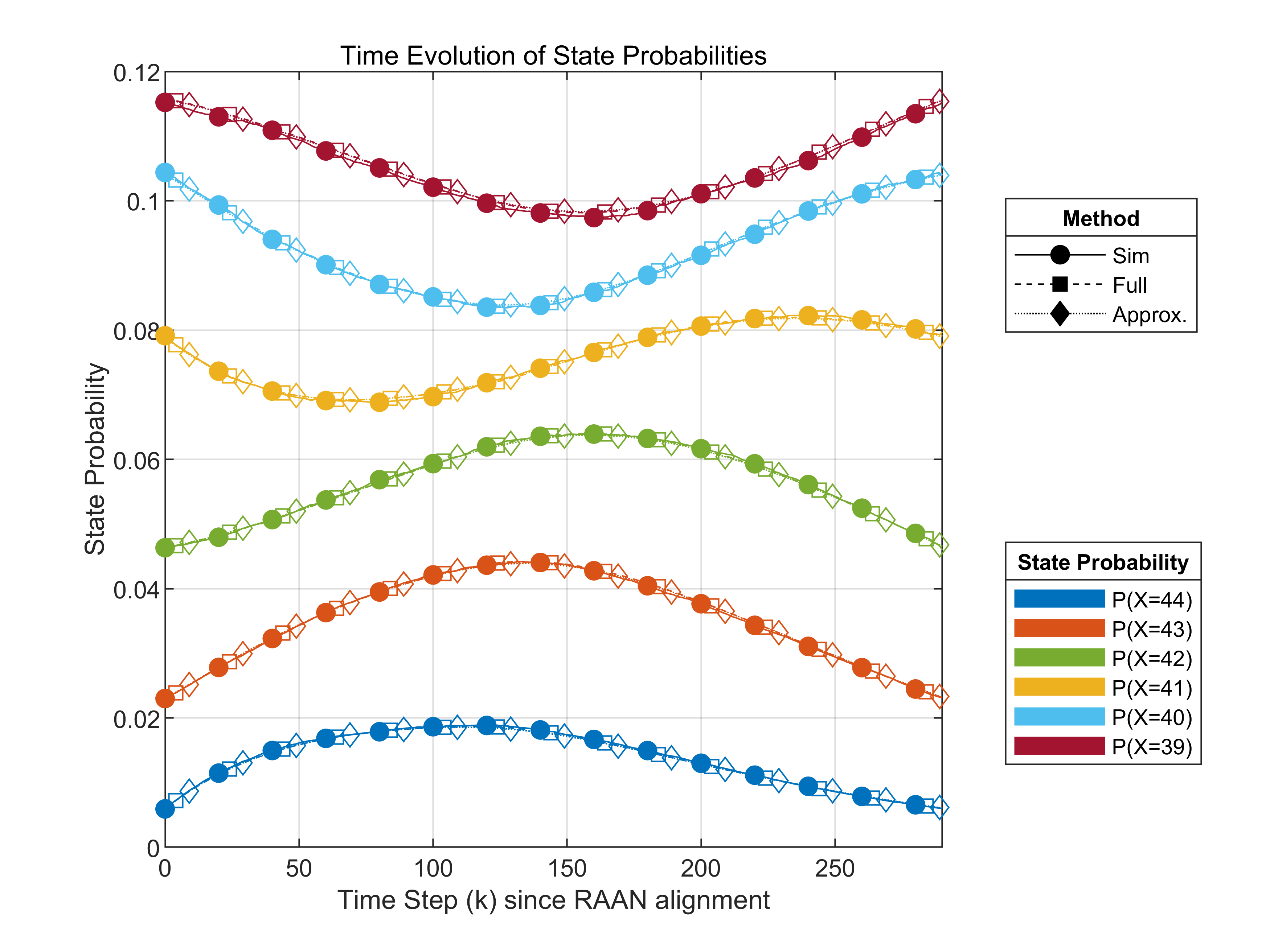}
    \caption{Time profile of selected in-plane state probabilities for the representative case}
    \label{fig:lead_time_test}
\end{figure}

The uniform lead-time assumption is also inaccurate. Reorder events are conditioned by RAAN alignment rather than occurring at random. Figure~\ref{fig:lead_time_test} shows the time evolution of selected in-plane state probabilities, $\boldsymbol{\pi}^\text{c}_k$, which follow non-trivial, $\tau_\text{c}$-periodic curves \textcolor{black}{rather than becoming independent of the RAAN-contact phase as a uniform lead-time approximation would imply}. This periodic behavior is consistent with the RAAN-conditioned lead-time behavior shown in our earlier indirect-strategy analysis~\cite{han2025indirect}, and indicates that replenishment is intrinsically coupled to orbital motion. Ignoring this effect overlooks the problem-specific multi-echelon dynamics of the system.

\subsubsection{Equivalence of Event-based and Time-Resolved Approaches}
As discussed in Sec.~\ref{sec:event_vs_cycle}, $(r,q,\tau)$ can be formulated in two ways: the event-based approach of \cite{han2025indirect} and the proposed time-resolved approach. Although they are expected to yield equivalent results, a formal mathematical proof is not tractable due to the coupled state variables, as shown in Eq.~\eqref{eq:stationary_parking}. To verify equivalence, we compared the solutions from both methods for the 100 test cases. The norm of the difference between $\pi^{\text{rc}_\text{p}}$ from the two approaches matches the numerical precision of solving $P \pi = \pi$ (about $10^{-14}$), confirming their equivalence.

\section{Optimization of the Spare Management Policy} \label{sec6}
\textcolor{black}{In this section, we apply the proposed framework to compare the hybrid, direct, and indirect strategies at the architecture level, and to examine how the optimal channel mix shifts with the operational environment.} Unless noted otherwise, all cost parameters and price models follow~\cite{han2025indirect} to enable direct comparison, and they are listed in Table~\ref{tab:opt_para}. Launch costs and lead-time parameters are taken from~\cite{SpaceInsider_cost,ainvest_rocketlab_2025,spaceflightnow_falcon9_2025}. The approximate analysis in Sec.~\ref{sec:approx_analysis} is used as the inner loop. We enforce the constraints in Eq.~\eqref{eq:opt_formulation} and solve the optimization with a genetic algorithm.

\begin{table}[hbt!]
\centering
\footnotesize
\setlength{\tabcolsep}{5pt}
\renewcommand{\arraystretch}{1.1}
\caption{Parameters for the optimization}
\label{tab:opt_para}
\begin{tabularx}{\linewidth}{@{}>{\raggedright\arraybackslash}X c c c@{}}
\hline \hline
Parameter & Notation & Value & Unit \\
\hline

Satellite manufacturing cost
  & $c_{\text{build}}$ 
  & $0.5$ 
  & M\$/satellite \\

In-orbit spares annual holding cost
  & $c_{\text{hold}_\text{c}}$ 
  & $0.5$ 
  & M\$/satellite/year \\

Parking spares annual holding cost
  & $c_{\text{hold}_\text{p}}$ 
  & $0.5$ 
  & M\$/satellite/year \\

Indirect Channel LV Lead Time Parameters
  & ($\tau_\text{lv,i}$,$\mu_\text{lv,i}$) 
  & ($20$, $20$)
  & days \\

Indirect Channel LV cost for full contract
  & $c_{\text{full,i}}$ 
  & $67$ 
  & M\$/launch \\

Indirect Channel Payload launch capacity 
  & $m_\text{payload,i}$ 
  & $18500$ 
  & kg \\

Direct Channel LV Lead Time Parameters
  & ($\tau_\text{lv,d}$,$\mu_\text{lv,d}$) 
  & ($10$, $10$)
  & days \\

Direct Channel LV cost for full contract
  & $c_{\text{full,d}}$ 
  & $7.5$ 
  & M\$/launch \\

Direct Channel Payload launch capacity
  & $m_\text{payload,d}$ 
  & $300$ 
  & kg \\  

Fuel cost per mass
  & $c_{\text{fuel}}$ 
  & $0.001$ 
  & M\$/kg/transfer \\

Non-fuel transfer cost
  & $c_{\text{trans}}$ 
  & $0.5$ 
  & M\$/transfer \\

Mass of satellite
  & $m_\text{sat}$ 
  & $150$ 
  & kg \\

Mass of transfer bus
  & $m_\text{bus}$ 
  & $100$ 
  & kg \\

Effective exhaust velocity 
  & $v_\text{ex}$ 
  & $2.16$ 
  & km/s \\

\hline
\hline
\end{tabularx}
\end{table}

\subsection{Baseline Scenario with Current Market Costs}
As a baseline, we consider a moderate annual failure rate ($\lambda_{\text{sat,yr}} = 0.05$) with a constellation resilience target of $S_c \leq 0.25$ and a parking availability target of $\mathbb{P}(X_\text{p} = 0) \leq 1/(N_{\text{sat}_\text{p}}+1)$. The optimization was performed using a genetic algorithm with cost parameters based on representative public estimates for Electron (direct channel) and Falcon 9 (indirect channel) launch vehicles. A summary of the optimal policies and their performance is presented in Table~\ref{tab:base_scn}, with a comparison to other strategies available in~\cite{han2025indirect}. \textcolor{black}{The primary result is that the optimal hybrid policy is functionally identical to the pure indirect policy: under the representative launch-cost assumptions, the channel-neutral framework selects the indirect channel as the preferred path, achieving the required resilience at less than half the total cost of the direct-only strategy.}


This outcome is driven by the significant cost-per-kilogram advantage of the Falcon 9. The optimization leverages this advantage by relying almost exclusively on the indirect channel; the direct-channel trigger is reached with negligible probability in the stationary distribution. The resulting strategy improves launch efficiency by using large indirect batches and a low parking-orbit altitude, which reduces launch frequency and shortens the RAAN-alignment interval. The build costs ($C_\text{build}$) are similar across all strategies, as the number of satellites produced is dictated by the same failure rate and resilience requirement. While the direct strategy satisfies the resilience constraint with a large margin ($S_c = 0.0591$), the optimal hybrid/indirect policy operates closer to the constraint boundary ($S_c = 0.2387$), \textcolor{black}{indicating a more cost-efficient operating point}.

\begin{table}[hbt!]
\caption{Summary of results for representative scenarios.}
\label{tab:base_scn}
\centering
\resizebox{\textwidth}{!}{%
\begin{tabular}{l c cccc cc}
\toprule
& & \multicolumn{4}{c}{Detailed Costs [M\$/day]} & \multicolumn{2}{c}{Constraints} \\
\cmidrule(lr){3-6} \cmidrule(lr){7-8}
Policy & $C_\text{total}$ [M\$/day] & $C_\text{build}$ & $C_\text{hold}$ & $C_\text{launch}$ & $C_\text{trans}$ & $S_\text{c}$ & $\mathbb{P}(X_\text{p} = 0)$ \\
\midrule
Direct & 0.9547 & 0.1094 & 0.0246 & 0.8207 & N/A & 0.0591 & N/A \\
Indirect & 0.4479 & 0.1082 & 0.1507 & 0.1575 & 0.0316 & 0.2387 & 0.0286 \\
Hybrid & 0.4479 & 0.1082 & 0.1507 & 0.1575 & 0.0316 & 0.2387 & 0.0286 \\
\bottomrule
\end{tabular}%
}
\end{table}

\subsection{Sensitivity to Direct Channel Launch Cost}
The baseline scenario confirmed the dominance of the indirect channel due to the cost-effectiveness of heavy-lift launches. This section explores the conditions under which the direct channel becomes the preferred option. To this end, we conduct a sensitivity analysis by varying the direct launch cost as $C_\text{full,d} \gets \beta C_\text{full,d}$, where the discount factor $\beta \in [0, 1]$.

\begin{figure}[h]
    \centering
    \includegraphics[width=.6\textwidth]{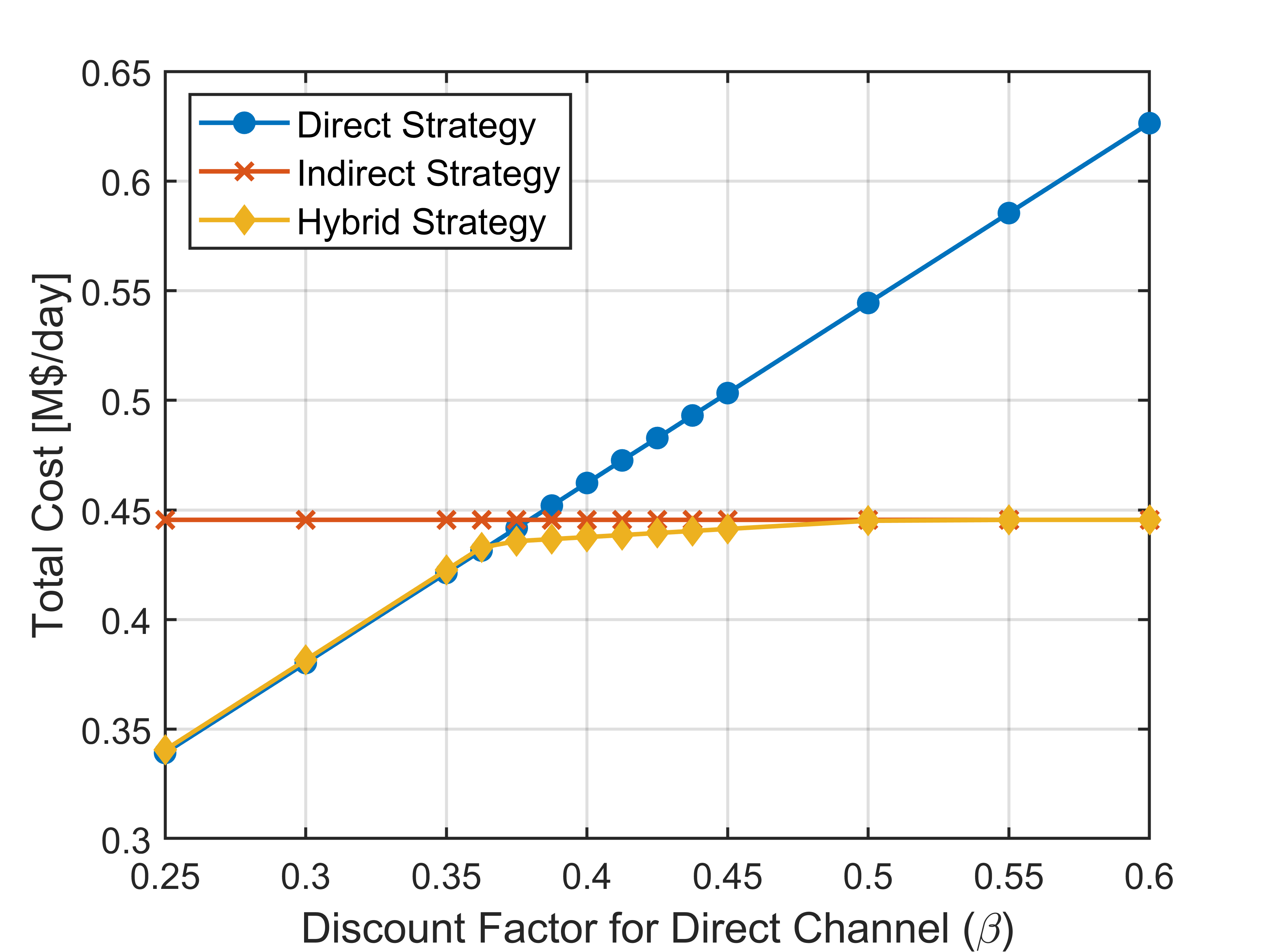}
    \caption{Total cost variation of each method at different $\beta$}
    \label{fig:cost_vs_beta}
\end{figure}

\begin{figure}[h]
    \centering
    \includegraphics[width=.6\textwidth]{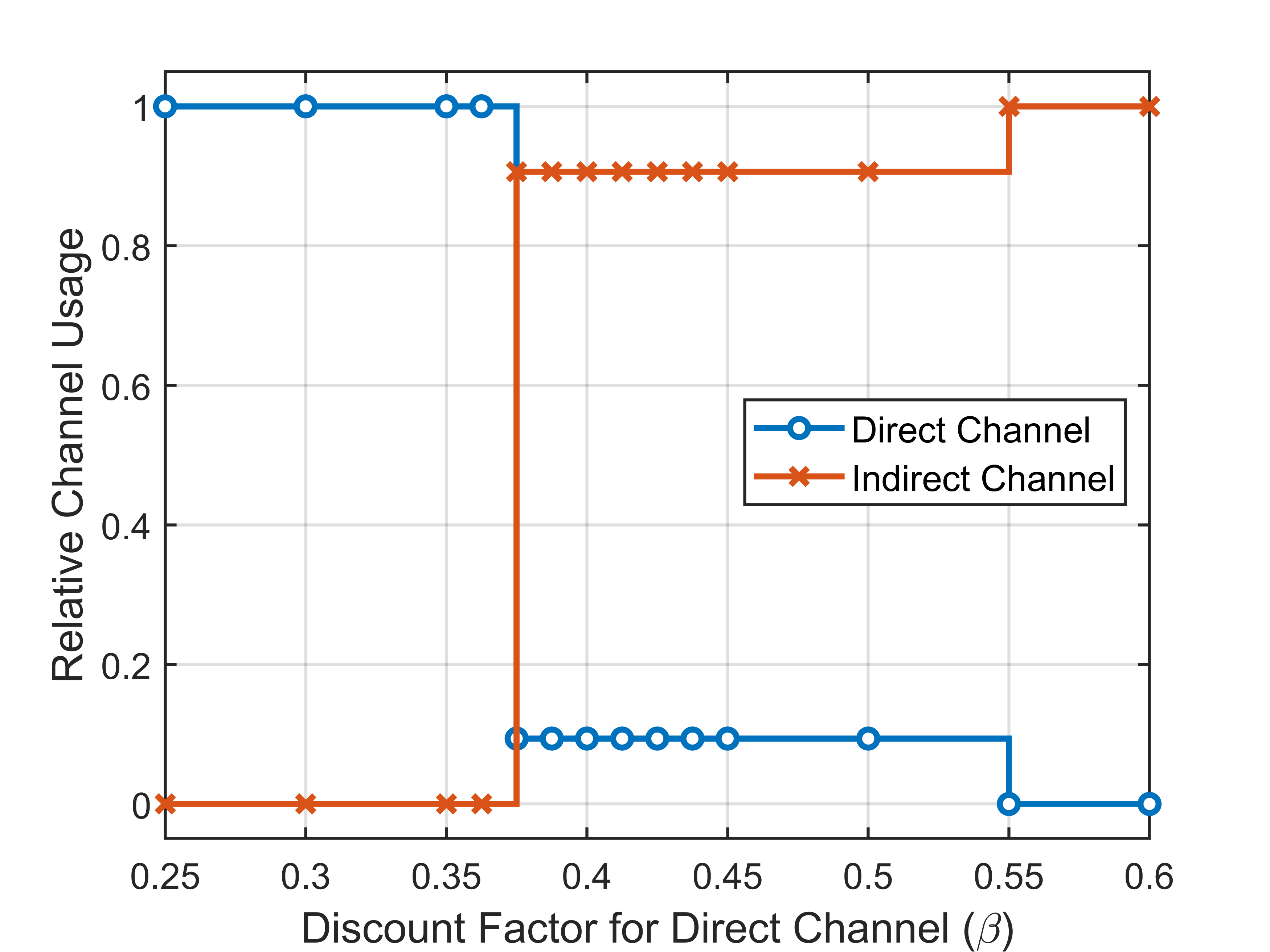}
    \caption{Relative channel usage in Hybrid Strategy at different $\beta$}
    \label{fig:usage_vs_beta}
\end{figure}

The results are shown in Figs.~\ref{fig:cost_vs_beta} and \ref{fig:usage_vs_beta}. The optimal hybrid strategy follows the lower cost envelope of the pure strategies, with a break-even point at $\beta^\ast \approx 0.38$. If the discount factor is less than $\beta^\ast$, the hybrid strategy mainly uses the direct channel, essentially becoming the direct strategy; the opposite is true for $\beta > \beta^\ast$. \textcolor{black}{In other words, the small launch vehicle would need to be discounted to roughly 38\% of its current price before the direct channel becomes economically preferred in this failure-rate regime, and across the entire $\beta$ range the channel-neutral framework consistently identifies the cheaper of the two pure strategies without any prior assignment of primary/secondary roles.}

Note that the cost of the hybrid strategy is slightly higher than the pure direct strategy even when direct channel usage is 100\%. This is because the model is constrained to maintain a minimal, but unused, parking orbit infrastructure ($N_{\text{orbit}_\text{p}} \geq 1$), which creates a small but unnecessary holding cost.

\subsection{Sensitivity to Satellite Failure Rate} \label{sec:sens_fail}
Satellite failure rates can vary significantly depending on system scale, with smaller satellites generally exhibiting higher failure rates than larger, more robust platforms \cite{2010_Gregory,langer2018reliability,jacklin2019small}. To assess the model's performance across this spectrum, we test a range of failure rates from 0.01 to 0.5 failures per satellite-year. For this analysis, the direct channel cost is held at the cost-neutral point identified previously ($\beta^\ast = 0.38$), isolating the impact of failure rate on the optimal strategy.

\begin{figure}[!h]
    \centering
    \includegraphics[width=.6\textwidth]{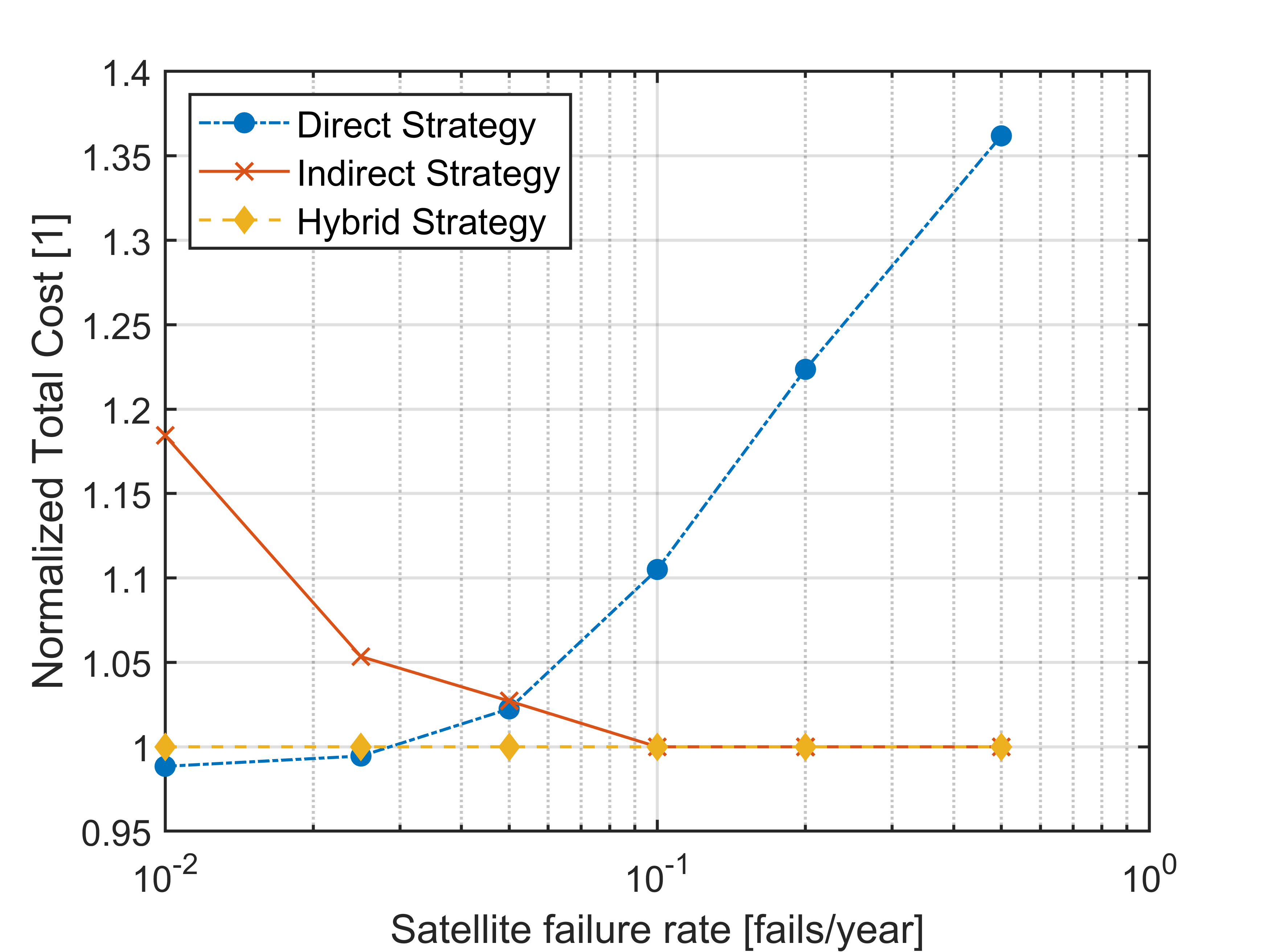}
    \caption{Total cost variation of each method at different $\lambda_{\text{sat,yr}}$}
    \label{fig:cost_vs_fail}
\end{figure}

Figure~\ref{fig:cost_vs_fail} shows the cost of each strategy normalized by the hybrid strategy's cost. The results confirm that the hybrid policy is the most cost-efficient across the entire range of failure rates. At very low failure rates, the demand for spares is very low as well. In this low-demand regime, the indirect channel is economically inefficient due to the high holding cost of large batches. The direct strategy is therefore superior, and the hybrid model correctly adopts this direct-dominant approach. Conversely, as the failure rate increases, the high volume of required launches makes the cost-per-kilogram of the indirect channel the dominant economic factor, and the optimal hybrid strategy converges to the indirect-dominant policy.  \textcolor{black}{This adaptive behavior is the key practical value of the framework: for a given launch-cost model and failure-rate assumption, it identifies whether a direct-dominant, indirect-dominant, or hybrid spare-management architecture is most suitable for early-phase design.}

\section{Conclusion} \label{sec7}
In this paper, we developed a time-resolved Markov chain framework for the early-phase analysis and design of a hybrid spare-management architecture that coordinates direct and indirect replenishment. The framework generalizes our earlier Markov chain analyses for the direct and indirect strategies~\cite{han2025direct,han2025indirect}. Constellation planes and parking orbits are modeled as $(r,q)$ and $(r,q,\tau)$ systems, and a coupled periodic steady state over the RAAN cycle is obtained by enforcing consistency via a fixed-point iteration. \textcolor{black}{Theoretically, this channel-neutral formulation captures the multi-echelon dynamics of dual replenishment paths without imposing aggregation assumptions.} From the periodic stationary distribution, we compute cost and resilience metrics, including total cost, expected shortage, and parking out-of-stock probability. We also provide an approximate analysis that preserves delay statistics while reducing the model size. Using this fast and accurate evaluation, we pose a cost minimization with resilience constraints, solve it with a genetic algorithm, and validate the results against Monte Carlo simulation. The optimization results show that the cost-effective indirect channel dominates under current market conditions. \textcolor{black}{From a practical perspective, the hybrid framework provides a quantitative early-phase assessment tool for examining how launch cost, response time, and replenishment throughput affect the suitability of a dual-channel spare-management architecture.}



\section*{Funding Sources}
This research was supported by the Advanced Technology R\&D Center at Mitsubishi Electric Corporation.

\section*{Declaration of Generative AI and AI-assisted technologies in the writing process}
During the preparation of this work, the authors used ChatGPT, Gemini, and Copilot to support initial background literature search, language polishing, and improving readability. After using these tools, the authors reviewed and edited the content as needed and take full responsibility for the content of the publication.

\bibliographystyle{elsarticle-num} 
\bibliography{references}

@article{ramasesh1991sole,
  title={Sole versus dual sourcing in stochastic lead-time (s, Q) inventory models},
  author={Ramasesh, Ranga V and Ord, J Keith and Hayya, Jack C and Pan, Andrew},
  journal={Management science},
  volume={37},
  number={4},
  pages={428--443},
  year={1991},
  publisher={INFORMS},
  doi={10.1287/mnsc.37.4.428}
}

@article{chiang1994sole,
  title={Sole sourcing versus dual sourcing under stochastic demands and lead times},
  author={Chiang, Chi and Benton, WC},
  journal={Naval Research Logistics (NRL)},
  volume={41},
  number={5},
  pages={609--624},
  year={1994},
  publisher={Wiley Online Library},
  doi={10.1002/1520-6750(199408)41:5<609::AID-NAV3220410503>3.0.CO;2-7}
}

@article{pan1991multiple,
  title={Multiple sourcing: the determination of lead times},
  author={Pan, Andrew C and Ramasesh, Ranga V and Hayya, Jack C and Ord, J Keith},
  journal={Operations Research Letters},
  volume={10},
  number={1},
  pages={1--7},
  year={1991},
  publisher={Elsevier},
  doi={10.1016/0167-6377(91)90079-5}
}

@book{schimpel2014dual,
    author       = {Schimpel, Ulrich},
    year         = {2010},
    title        = {Dual sourcing : with arbitrary stochastic demand and stochastic lead times},
    doi          = {10.5445/KSP/1000018435},
    publisher    = {{KIT Scientific Publishing}}
}

@article{SUNG_deploy,
title = {Optimal deployment of satellite mega-constellation},
journal = {Acta Astronautica},
volume = {202},
pages = {653-669},
year = {2023},
issn = {0094-5765},
doi = {10.1016/j.actaastro.2022.10.027},
author  = {Sung, Taehyun and Ahn, Jaemyung}
}

@article{Kim2025_dual,
author = {Kim, Jaewoo and Ahn, Jaemyung and Sung, Taehyun},
title = {Replenishment Strategy for Satellite Constellation with Dual Supply Modes},
journal = {Journal of Spacecraft and Rockets},
volume = {62},
number = {5},
pages = {1567-1583},
year = {2025},
doi = {10.2514/1.A36281}
}

@misc{starlink1,
  author       = {Jonathan McDowell},
  title        = {Starlink Statistics},
  howpublished = {\url{https://planet4589.org/space/con/star/stats.html}},
  year         = {2025},
  note         = {Accessed: 2025-07-14}
}

@misc{starlink2,
  author       = {{Space.com}},
  title        = {Starlink satellites: Facts, tracking and impact on astronomy},
  howpublished = {\url{https://www.space.com/spacex-starlink-satellites.html}},
  year         = {2025},
  note         = {Accessed: 2025-08-07}
}

@misc{kuiper2025,
  author       = {{CNBC}},
  title        = {Amazon launches second batch of Kuiper internet satellites, taking on Elon Musk's Starlink},
  howpublished = {\url{https://www.cnbc.com/2025/06/23/amazon-kuiper-satellites-musk-starlink.html}},
  year         = {2025},
  note         = {Accessed: 2025-07-14}
}

@misc{oneweb,
  author       = {{Space.com}},
  title        = {Arianespace to launch new fleet of OneWeb internet satellites tonight. Here's how to watch},
  howpublished = {\url{https://www.space.com/arianespace-oneweb-5-satellites-launch-webcast}},
  year         = {2021},
  note         = {Accessed: 2025-07-14}
}

@misc{guowang2025,
  author       = {{Telecoms Tech News}},
  title        = {China launches first satellites for GuoWang constellation},
  howpublished = {\url{https://www.telecomstechnews.com/news/china-launches-first-satellites-guowang-constellation/}},
  year         = {2024},
  note         = {Accessed: 2025-07-14}
}

@inproceedings{Pachler2021Comp,
  author    = {Pachler, Nils and del Portillo, Inigo and Crawley, Edward F. and Cameron, Bruce G.},
  title     = {An Updated Comparison of Four Low Earth Orbit Satellite Constellation Systems to Provide Global Broadband},
  booktitle = {2021 IEEE International Conference on Communications Workshops (ICC Workshops)},
  year      = {2021},
  pages     = {1--7},
  doi       = {10.1109/ICCWorkshops50388.2021.9473799}
}

@article{Pachler2024Comp,
  author  = {Pachler, Nils and Crawley, Edward F. and Cameron, Bruce G.},
  title   = {Flooding the Market: Comparing the Performance of Nine Broadband Megaconstellations},
  journal = {IEEE Wireless Communications Letters},
  year    = {2024},
  volume  = {13},
  number  = {9},
  pages   = {2397--2401},
  doi     = {10.1109/LWC.2024.3416531}
}

@article{selva2017distributed,
author={Selva, Daniel and Golkar, Alessandro and Korobova, Olga and Cruz, Ignasi Lluch i and Collopy, Paul and de Weck, Olivier L.},
year={2017},
title={Distributed Earth Satellite Systems: What Is Needed to Move Forward?},
journal={Journal of Aerospace Information Systems},
volume={14},
number={8},
pages={412--438},
doi={10.2514/1.I010497},
}

@misc{potter2023mooreslaw,
  author       = {{American Enterprise Institute}},
  title        = {Moore's Law Meets Musk's Law: The Underappreciated Story of SpaceX and the Stunning Decline in Launch Costs},
  howpublished = {\url{https://www.aei.org/articles/moores-law-meet-musks-law-the-underappreciated-story-of-spacex-and-the-stunning-decline-in-launch-costs/}},
  year         = {2024},
  note         = {Accessed: 2025-07-14}
}

@article{BOUWMEESTER2022108288,
title = {Improving CubeSat reliability: Subsystem redundancy or improved testing?},
journal = {Reliability Engineering \& System Safety},
volume = {220},
pages = {108288},
year = {2022},
issn = {0951-8320},
doi = {10.1016/j.ress.2021.108288},
author = {J. Bouwmeester and A. Menicucci and E.K.A. Gill},
}

@article{luu2022orbit,
  title={On-orbit servicing system architectures for proliferated low-earth-orbit constellations},
  author={Luu, Michael A and Hastings, Daniel E},
  journal={Journal of Spacecraft and Rockets},
  volume={59},
  number={6},
  pages={1946--1965},
  year={2022},
  publisher={American Institute of Aeronautics and Astronautics},
  doi={10.2514/1.A35393}
}

@article{jakob2019optimal,
  title={Optimal satellite constellation spare strategy using multi-echelon inventory control},
  author={Jakob, Pauline and Shimizu, Seiichi and Yoshikawa, Shoji and Ho, Koki},
  journal={Journal of Spacecraft and Rockets},
  volume={56},
  number={5},
  pages={1449--1461},
  year={2019},
  publisher={American Institute of Aeronautics and Astronautics},
  doi={10.2514/1.A34387}
}

@mastersthesis{sumter2003optimal,
  author       = {Bradley R. Sumter},
  title        = {Optimal Replacement Policies for Satellite Constellations},
  school       = {Air Force Institute of Technology},
  department   = {Department of Operational Sciences},
  year         = {2003}
}

@inproceedings{cornara1999satellite,
  title={Satellite constellation launch, deployment, replacement and end-of-life strategies},
  author={Cornara, Stefania and Beech, Theresa and Bell{\'o}-Mora, Miguel and Martinez de Aragon, Antonio},
  booktitle={13th Annual AIAA/USU Conference on Small Satellites},
  year={1999}
}

@article{1966dishon,
 ISSN = {00401706},
 author = {Menachem Dishon and George H. Weiss},
 journal = {Technometrics},
 number = {3},
 pages = {399--409},
 publisher = {[Taylor \& Francis, Ltd., American Statistical Association, American Society for Quality]},
 title = {A Communications Satellite Replenishment Policy},
 urldate = {2025-07-14},
 volume = {8},
 year = {1966},
 doi = {10.1080/00401706.1966.10490373}
}

@book{hillier2015introduction,
  title={Introduction to operations research},
  author={Hillier, Frederick S and Lieberman, Gerald J},
  year={2015},
  publisher={McGraw-Hill}
}

@article{walker1984satellite,
  title={Satellite constellations},
  author={Walker, John G},
  journal={Journal of the British Interplanetary Society},
  volume={37},
  pages={559--572},
  year={1984}
}

@book{prussing1993orbital,
  title={Orbital mechanics},
  author={Prussing, John E and Conway, Bruce A},
  year={1993},
  publisher={Oxford University Press, USA}
}

@inproceedings{han2024analysis,
  author    = {Han, Seungyeop and Noro, Takumi and Ho, Koki},
  title     = {Analysis and Design of Satellite Constellation Spare Strategy Using Markov Chain},
  booktitle = {2024 AAS/AIAA Astrodynamics Specialist Conference},
  year      = {2024},
  note      = {Broomfield, Colorado, USA, August 2024}
}

@article{han2025direct,
  title        = {Analysis and Design of Spare Strategy for Large-Scale Satellite Constellation Using Direct Insertion under $(r,q)$ policy},
  author       = {Han, Seungyeop and Grieser, Zachary and Yoshikawa, Shoji and Noro, Takumi and Suda, Takumi and Ho, Koki},
  year         = {2025},
  journal      = {arXiv preprint arXiv:2509.10585},
  doi          = {10.48550/arXiv.2509.10585}
}

@article{han2025indirect,
  title   = {Spare Strategy Analysis and Design for Mega Satellite Constellations Using {Markov} Chain},
  author  = {Han, Seungyeop and Grieser, Zachary and Yoshikawa, Shoji and Noro, Takumi and Suda, Takumi and Ho, Koki},
  journal = {Journal of Spacecraft and Rockets},
  year    = {2026},
  note    = {, accepted for publication},
  eprint  = {2509.09957},
  archivePrefix = {arXiv}
}

@INPROCEEDINGS{Ron_resilience,
  author={Burch, Ron},
  booktitle={MILCOM 2013 - 2013 IEEE Military Communications Conference}, 
  title={A Method for Calculation of the Resilience of a Space System}, 
  year={2013},
  volume={},
  number={},
  pages={1002-1007},
  doi={10.1109/MILCOM.2013.174}}

@misc{dod_resilience_2011,
  author       = {{U.S. Department of Defense}},
  title        = {Fact Sheet: Resilience of Space Capabilities},
  year         = {2011},
  howpublished = {\url{ https://web.archive.org/web/20210323071718/https://archive.defense.gov/home/features/2011/0111_nsss/docs/DoD\%20Fact\%20Sheet\%20-\%20Resilience.pdf }},
  note         = {Accessed: 2026-06-18}
}

@misc{SpaceInsider_cost,
  author       = {{SpaceInsider}},
  title        = {How Much Does It Cost to Launch a Rocket? [By Type \& Size]},
  year         = {2025},
  howpublished = {\url{https://spaceinsider.tech/2023/08/16/how-much-does-it-cost-to-launch-a-rocket/}},
  note         = {Accessed: 2025-07-11}
}

@misc{ainvest_rocketlab_2025,
  author       = {{Ainvest}},
  title        = {Rocket Lab's 48-Hour Launch Turnaround: A Strategic Play for Dominance in the SmallSat Launch Market},
  year         = {2025},
  howpublished = {\url{https://www.ainvest.com/news/rocket-lab-48-hour-launch-turnaround-strategic-play-dominance-smallsat-launch-market-2506/}},
  note         = {Accessed: 2025-07-11}
}

@misc{spaceflightnow_falcon9_2025,
  author       = {{Spaceflight Now}},
  title        = {SpaceX breaks launchpad turnaround record with midnight Starlink flight},
  year         = {2025},
  howpublished = {\url{ https://spaceflightnow.com/2025/06/27/live-coverage-spacex-to-launch-27-starlink-satellites-on-falcon-9-rocket-from-cape-canaveral-4/ }},
  note         = {Accessed: 2025-07-11}
}

@article{2010_Gregory,
  title   = {Statistical reliability analysis of satellites by mass category: Does spacecraft size matter?},
  author  = {Dubos, Gregory F. and Castet, Jean-Francois and Saleh, Joseph H.},
  journal = {Acta Astronautica},
  volume  = {67},
  number  = {5},
  pages   = {584--595},
  year    = {2010},
  issn    = {0094-5765},
  doi     = {10.1016/j.actaastro.2010.04.017}
}

@phdthesis{langer2018reliability,
  title={Reliability assessment and reliability prediction of CubeSats through system level testing and reliability growth modelling},
  author={Langer, Martin},
  year={2018},
  school={Technische Universit{\"a}t M{\"u}nchen}
}

@techreport{jacklin2019small,
  author       = {Jacklin, Stephen A.},
  title        = {Small‑Satellite Mission Failure Rates},
  institution  = {NASA Ames Research Center},
  type         = {NASA Technical Memorandum},
  number       = {NASA/TM‑2018‑220034},
  address      = {Moffett Field, CA},
  year         = {2019},
  month        = {March}
}

@book{greenbaum1997iterative,
  title     = {Iterative Methods for Solving Linear Systems},
  author    = {Greenbaum, Anne},
  year      = {1997},
  publisher = {SIAM},
  doi       = {10.1137/1.9781611970937}
}
\end{document}